\documentclass[12pt]{amsart}
\usepackage{amsfonts, amsmath, amsthm}
\usepackage{epsfig}
\usepackage{psfrag}

\newtheorem{pro}{Proposition}[section]
\newtheorem{thm}[pro]{Theorem}
\newtheorem{lem}[pro]{Lemma}
\newtheorem{clm}[pro]{Claim}
\newtheorem{example}[pro]{Example}

\newtheorem{cor}[pro]{Corollary}

\theoremstyle{definition}
\newtheorem{dfn}[pro]{Definition}

\theoremstyle{remark}

\title{Heegaard Splittings with Boundary and Almost Normal Surfaces}
\date{\today}
\address{Mathematics Department, California Polytechnic State University, San Luis Obispo}
\email{dbachman@calpoly.edu}
\author{David Bachman}

\begin{document}
\begin{abstract}
This paper generalizes the definition of a Heegaard splitting to unify the concepts of thin position for 3-manifolds \cite{st:94}, thin position for knots \cite{gabai:87}, and normal and almost normal surface theory \cite{haken:61}, \cite{rubinstein:93}. This gives generalizations of theorems of Scharlemann, Thompson, Rubinstein, and Stocking. In the final section, we use this machinery to produce an algorithm to determine the bridge number of a knot, provided thin position for the knot coincides with bridge position. We also present several algorithmic and finiteness results about Dehn fillings with small Heegaard genus.
\end{abstract}
\maketitle

\pagestyle{myheadings}
\markboth{DAVID BACHMAN}{HEEGAARD SPLITTINGS WITH BOUNDARY [Top. and Appl. 116(2001) 153-184]}


\section{Introduction.}
\label{s:introduction}

A Morse function on a closed 3-manifold, $M$, is a generic height function, $h:M \rightarrow I$. The set, $\mathcal F$, of level sets of such a function forms a ``singular" foliation of $M$. Generic elements of $\mathcal F$ are closed surfaces. A component of a non-generic element is either a point, a pair of surfaces which meet at a point, or one surface which touches itself at a point. In the next section we will define a complexity, $c$, on generic elements of $\mathcal F$, which has the property that if $F_2 \in \mathcal F$ is obtained from $F_1 \in \mathcal F$ by a compression, then $c(F_2)<c(F_1)$. This allows us to talk about ``maximal" and ``minimal" leaves, which are simply elements of $\mathcal F$ which correspond to local maxima and minima of $c$. We then see that the submanifolds of $M$ between maximal and minimal leaves are standard objects of 3-manifold topology, called {\it compression bodies}. If two adjacent compression bodies, $W$ and $W'$, have a maximal leaf, $F$, between them then $W \cup _F W'$ is called a {\it Heegaard splitting}, and $F$ is referred to as a {\it Heegaard surface}. From this, we immediately deduce that the minimal leaves of $\mathcal F$ break up $M$ into submanifolds, for which the maximal leaves are Heegaard surfaces. This is precisely the picture of a 3-manifold presented by Scharlemann and Thompson in \cite{st:94}.

In Section \ref{s:BoundaryHS} we generalize this picture to manifolds with boundary. That is, suppose now that $h:M \rightarrow I$ is a Morse function on a compact 3-manifold with nonempty boundary, such that $h|_{\partial M}$ is also a Morse function. Once again, we denote the set of level sets of $h$ as $\mathcal F$. Now, a generic element of $\mathcal F$ is generally a surface with non-empty boundary. After a slight modification of our complexity, $c$, we can once again talk about maximal and minimal leaves of $\mathcal F$. But the submanifolds of $M$ between maximal and minimal leaves are no longer compression bodies, in the usual sense. This motivates us to define a {\it $\partial$-compression body} to be just such a submanifold. And if adjacent $\partial$-compression bodies, $W$ and $W'$, share a maximal leaf, $F$, we refer to $W \cup _F W'$ as a {\it $\partial$-Heegaard splitting}. Our picture is now exactly the same as before: the minimal leaves of $\mathcal F$ break up $M$ into submanifolds, for which the maximal leaves are $\partial$-Heegaard surfaces.

Sections \ref{s:Heegaard} and \ref{s:BoundaryHS} also present various notions of nontriviality for Heegaard and $\partial$-Heegaard splittings, namely the concepts of {\it strong irreducibility}, and, somewhat weaker, {\it quasi-strong irreducibility}. As these definitions can be quite difficult to get a feel for, we present several illustrative examples in Section \ref{ex}. It would be well worth the reader's time to get a good understanding of each example presented. Some of the main theorems presented later in the paper are simply generalizations of these examples to arbitrary manifolds.

Section \ref{s:closed} begins by presenting a complexity for height functions on $M$. We then show that if $\mathcal F$ is the set of level sets of a height function which minimizes this measure of complexity, then the maximal leaves of $\mathcal F$ are strongly irreducible ($\partial$-)Heegaard splittings for the submanifolds between the minimal leaves. From this, it follows that the minimal leaves are incompressible and $\partial$-incompressible in $M$. If $M$ is closed, these are the results of Scharlemann and Thompson from \cite{st:94}.

Next, we turn to embedded 1-manifolds in $M$. Section \ref{s:knots} examines the following question: What happens when we begin with a minimal height function, $h$, and isotope some 1-manifold, $K$, so that the complexity of $h|_{M-N(K)}$ (where $N(K)$ denotes a small neighborhood of $K$) is as small as possible? If $K$ is in such a position, we call $K$ {\it mini-Lmax}, which is a generalization of the minimax complexity. If $M$ is homeomorphic to $S^3$, then mini-Lmax is very similar to the {\it thin position} of Gabai \cite{gabai:87}. In \cite{thompson:97}, Thompson proves the following theorem:

\begin{thm}
\label{thompson}
{\rm [Thompson]} Suppose $K$ is a knot in $S^3$. Then thin position for $K$ is bridge position, or there is a meridional, incompressible, planar surface in the complement of $K$.
\end{thm}

If $K$ is the unknot, this theorem is trivially true. However, in some sense, the unknot lacks some of the nice properties of a knot in thin position. It therefore does no harm to rule out the unknot from the statement of Theorem \ref{thompson}, and there are some aesthetic reasons for doing this. One (albeit overly complicated) way to do this is to define a strongly irreducible Heegaard surface for $S^3$ to be any embedded 2-sphere, and restate Theorem \ref{thompson} as:

\medskip
\noindent
{\bf Theorem \ref{thompson}'} [Thompson] {\it Suppose $K$ is a knot in $S^3$, and $H$ is a strongly irreducible Heegaard surface. If $K$ cannot be isotoped onto $H$, then thin position for $K$ is bridge position, or there is a meridional, incompressible, $\partial$-incompressible, planar surface in the complement of $K$.}
\medskip

Section \ref{s:knots} ends with the following generalization of the above theorem:

\medskip
\noindent
{\bf Theorem \ref{t:generalthompson}} {\it Suppose $K$ is a knot in a closed, orientable, irreducible 3-manifold, $M$, and $H$ is a strongly irreducible Heegaard surface. If $K$ cannot be isotoped onto $H$, then mini-Lmax position for $K$ is bridge position, or there is a meridional, incompressible, $\partial$-incompressible surface in the complement of $K$, which has genus less than or equal to that of $H$.}
\medskip

The remainder of the paper deals with relating the above results to the theory of normal surfaces. After a brief review of this theory in Section \ref{s:normal}, we turn to the following question in Section \ref{s:nlan}: If we begin with a 1-vertex psuedo-triangulation of $M$, and make the 1-skeleton mini-Lmax away from the vertex, then what to the maximal and minimal leaves of $\mathcal F$ look like inside each tetrahedron? After a careful analysis, we find that the minimal leaves are a union of triangles and quadralaterals (i.e. a normal surface). We also find that the maximal leaves are a union of triangles and quadralaterals, except for exactly one exceptional piece (i.e an almost normal surface), where all possible exceptional pieces can easily be classified. The section ends with a generalization of a Theorem of Rubinstein \cite{rubinstein:93} and Stocking \cite{stocking:96}, that any strongly irreducible ($\partial$-)Heegaard splitting can be isotoped to be almost normal.

The last section (\ref{s:app}) focuses on applications of the above results. In most of this section, the manifolds which we consider are the complements of knots in arbitrary 3-manifolds. After triangulating in a special way, and making the 1-skeleton mini-Lmax, we discover the existence of many interesting normal and almost normal surfaces. If we are in the special case of a hyperbolic knot for which thin position is the same as bridge position, then this gives an algorithm to determine bridge number of that knot.

We also show that our existence results, when combined with a recent finiteness result of Jaco and Sedgwick \cite{js:98}, give several interesting results about Dehn filling. Suppose $X$ is a compact, irreducible, orientable 3-manifold with a single boundary component, homeomorphic to a torus. A {\it Dehn filling} of $X$ refers to the process of constructing a new manifold, by identifying the boundary of $X$ with the boundary of a solid torus. One of our more surprising results is the following:

\medskip
\noindent
{\bf Corollary \ref{Dehn}} {\it For all but finitely many Dehn fillings of $X$, the core of the attached solid torus can be isotoped onto every strongly irreducible Heegaard surface.}
\medskip

We are also able to reproduce some of the algorithmic results of Jaco and Sedgwick from \cite{js:98}. In particular, we give new algorithms to determine if $X$ is the complement of a knot in $S^3$, a lens space, or $S^2 \times S^1$.

\section{Heegaard Splittings and Morse Theory.}
\label{s:Heegaard}

In this section we review some of the basic definitions and facts about Heegaard splittings, and review their relationship to Morse theory.

$M$ will always denote a compact, orientable 3-manifold. An embedded 2-sphere in $M$ is {\it essential} if it does not bound a 3-ball. A manifold which does not contain an essential 2-sphere is {\it irreducible}. It will be assumed that all 3-manifolds considered in this paper are irreducible.

Let $F$ denote a compact, orientable surface, embedded in $M$ (possibly, $F \subset \partial M$). An {\it essential curve} on $F$ is an embedded loop, which does not bound a disk on $F$. A {\it compressing disk} for $F$ is a disk, $D$, embedded in $M$, such that $D \cap F= \partial D$, and $\partial D$ is essential on $F$. If such a disk exists, then $F$ is {\it compressible}; otherwise, it is {\it incompressible}.

Now, suppose $D$ is a compressing disk for $F$. Then there exists an embedding, $\phi :D^2 \times I \rightarrow M$, such that $D=\phi (D^2 \times \{1/2\})$, and $F \cap  \phi (D^2 \times I)= \phi (\partial D^2 \times I)$. {\it Surgery of $F$ along D} simply refers to the process of removing $\phi(\partial D^2 \times I)$ from $F$, and replacing it with $\phi(D^2 \times \partial I)$. We shall also refer to a surgery of $F$ as a {\it compression} of $F$.

Let $h:M \rightarrow [0,1]$ be a Morse function, where we require that $\partial M \subset h^{-1}(0) \cup h^{-1}(1)$ (if $\partial M \ne \emptyset$). $h$ determines a singular foliation, $\mathcal F$, of $M$ in the usual way, where the leaves of $\mathcal F$ are the inverse images of points in $[0,1]$, and a generic leaf is a compact, embedded surface. For each $t \in [0,1]$, let $\mathcal F_t =h^{-1}(t)$. We now define a complexity on $\mathcal F_t$, assuming $t$ is not a critical value of $h$.

Suppose $\mathcal F_t^i$ is a component of $\mathcal F_t$. Define $c(\mathcal F_t^i)$ to be $0$ if $\mathcal F_t^i$ is a sphere, and $1-\chi (\mathcal F_t^i)$ otherwise. Let $c(\mathcal F_t)=\sum _i c(\mathcal F_t^i)$, where the sum is taken over all components of $\mathcal F_t$. This measure of complexity will decrease if we see any compression of $\mathcal F_t$, and it will be $0$ if and only if $\mathcal F_t$ is a collection of spheres.

Let $\{ s_i \}$ be some collection of points in [0,1], such that there is exactly one element of this set between any two consecutive critical values of $h$. It is important to note that we can obtain $\mathcal F_{s_i}$ from $\mathcal F_{s_{i-1}}$ by either adding or removing a 2-sphere, or by compressing or ``de-compressing" (the reverse of a compression). Hence, we can build $M$ by a handle decomposition, where the surface $\mathcal F_{s_i}$ is the boundary of the manifold we get after adding the $i$th handle.

Now, let $\{ t_i \}$ be some subcollection of $\{ s_i \}$ such that $\mathcal F_{t_i}$ differs from $\mathcal F_{t_{i+1}}$ by exactly one compression or de-compression (and possibly several 2-sphere components). So, by definition we have $c(\mathcal F_{t_i}) \ne c(\mathcal F_{t_{i+1}})$. We say that a {\it local maximum occurs at} $t_i$ if $c(\mathcal F_{t_i})>c(\mathcal F_{t_{i-1}})$ and $c(\mathcal F_{t_i})>c(\mathcal F_{t_{i+1}})$. We can define a {\it local minimum} in a similar manner. If a local maximum (minimum) occurs at $t_i$, then we refer to $\mathcal F_{t_i}$ as a {\it maximal (minimal) leaf} of $\mathcal F$.

We now ask the following question: What do submanifolds of $M$ between consecutive maximal and minimal leaves look like? This is a standard object of 3-manifold topology, called a {\it compression body}, which we shall define in several ways.

We say a separating surface, $F$, is {\it completely compressible to one side} if there exists a collection of disjoint compressing disks for $F$ on one side, such that surgery along every disk in this collection yields a collection of spheres which bound balls, or yields a surface which is parallel to some subsurface of $\partial M$.

A {\it compression body} is a 3-manifold $W$, such that $\partial W$ is the union of 2 subsurfaces, denoted $\partial _+ W, \partial _- W$, such that $\partial _+W$ is completely compressible to one side, and when compressed, is parallel to $\partial _-W$ (if $\partial _- W \ne \emptyset$), or is a collection of 2-spheres which bound balls (if $\partial _-W =\emptyset$). We also insist that all compression bodies are nontrivial, in the sense that $\partial _+ W$ is not homeomorphic to $\partial _- W$. In other words, we are not allowing a compression body to be a product.

Another description of a compression body is any 3-manifold that can be built up in the following way: Begin with a closed, orientable surface, $F$, and form the product $F \times I$. Denote $F \times \{ 0 \}$ by $\partial _+ W$. Now, add a non-empty collection of 2-handles to $F \times \{ 1 \}$, and cap off any resulting 2-sphere boundary components by 3-balls. We denote $\partial W \backslash \partial _+ W$ by $\partial _- W$. It follows that $\partial _- W$ is incompressible in $W$.

A {\it Heegaard splitting} of a manifold, $M$, is a decomposition into two compression bodies, $W$ and $W'$, such that $W \cap W'=\partial _+ W=\partial _+ W'=F$. We denote such a splitting as $W \cup _F W'$. Another way to say this is that there is a surface, $F \subset M$ which is completely compressible to both sides. It is easy to show that every 3-manifold posesses infinitely many Heegaard splittings. In 1987, Casson and Gordon \cite{cg:87} introduced a notion of non-triviality for Heegaard splittings. A {\it strongly irreducible Heegaard splitting} is one which has the property that every compressing disk for $F$ in $W$ must have non-empty intersection with every compressing disk for $F$ in $W'$.

One of the main theorems that makes strongly irreducible Heegaard splittings useful is the following:

\begin{thm}
\label{t:cg}
If $W \cup _F W'$ is a strongly irreducible splitting of $M$, then $\partial M$ is incompressible in $M$.
\end{thm}

This was originally proved by Casson and Gordon in \cite{cg:87}, by using a Lemma of Haken (Lemma 1.1 in the Casson-Gordon paper. See \cite{haken:68} for the original Lemma). In the next section we generalize the concept of a Heegaard Splitting, and present an analogous result.
Our proof (in Appendix \ref{a:bcg}) will not use the Haken Lemma, and is general enough to include a new, simpler proof of Theorem \ref{t:cg}.

Let us now go back to the Morse function, $h:M \rightarrow I$, and the singular foliation, $\mathcal F$, which it defines. As we move from a maximal to a minimal leaf of $\mathcal F$ we see a sequence of compressions, and 2-spheres being capped off. Hence, a region between consecutive maximal and minimal leaves is precisely a compression body. The minimal leaves therefore break $M$ up into submanifolds, where each such submanifold contains a single maximal leaf, which is a Heegaard splitting. This is the point of view presented by Scharlemann and Thompson in \cite{st:94}.

\section{$\partial$-Heegaard Splittings.}
\label{s:BoundaryHS}

We now ask the question: What happens if we have a Morse function, $h:M \rightarrow [0,1]$, where $\partial M$ is not contained in $h^{-1}(0) \cup h^{-1}(1)$? In particular, what happens when $h$ restricted to $\partial M$ is a Morse function? In this case, a generic leaf, $h^{-1}(t)$, is not necessarily a closed surface. And as $t$ changes, we may see $h^{-1}(t)$ change in ways other than compression, de-compression, and addition and subtraction of 2-spheres.

\begin{figure}[htbp]
    \psfrag{F}{$F$}
    \psfrag{M}{$\partial M$}
    \psfrag{a}{(a)}
    \psfrag{b}{(b)}
    \psfrag{c}{(c)}
        \vspace{0 in}
        \begin{center}
        \epsfxsize=4 in
        \epsfbox{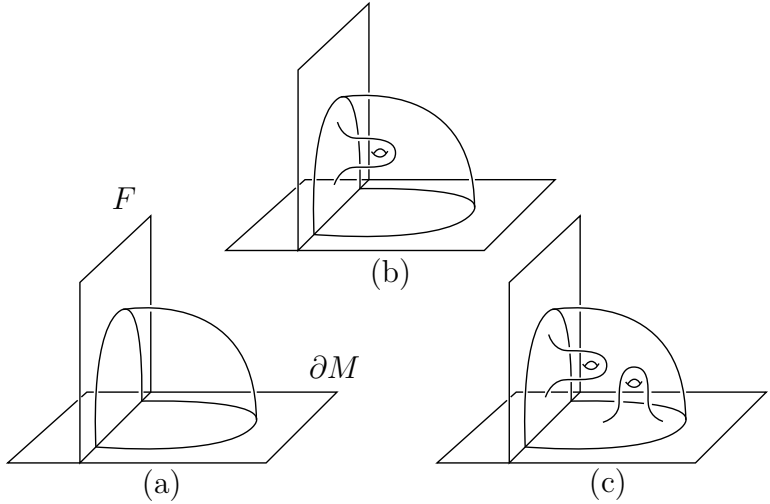}
        \caption{(a) Not a $\partial$-compression. (b) A $\partial$-compression which is not honest. (c) An honest $\partial$-compression.}
        \label{f:bcomp}
        \end{center}
        \end{figure}

To completely describe what may happen, we must first generalize the definitions given in the previous section. Suppose $(F,\partial F) \subset (M,\partial M)$. An {\it essential arc} $(\alpha, \partial \alpha) \subset (F, \partial F)$ is an embedded arc, such that there is no arc, $\beta \subset \partial F$, where $\alpha \cup \beta$ bounds a disk on $F$. A $\partial${\it -compressing disk} for $F$ is a disk, $D$, embedded in $M$, such that $\partial D=\alpha \cup \beta$, $D \cap F= \alpha$, $D \cap \partial M =\beta$, and $\alpha$ is an essential arc on $F$. If such a disk exists, then $F$ is $\partial${\it -compressible}; otherwise, it is $\partial${\it -incompressible}. If, in addition, $\beta$ is essential in $\partial M \backslash \partial F$, then we say $D$ is an {\it honest} $\partial$-compressing disk (see figure \ref{f:bcomp}).

Now, suppose $D$ is a $\partial$-compressing disk for $F$. Then there exists an embedding, $\phi :D^2 \times I \rightarrow M$, such that $D=\phi (D^2 \times \{1/2\})$, $F \cap  \phi (D^2 \times I)= \phi (\alpha \times I)$, and $\partial M \cap \phi (D^2 \times I)=\phi (\beta \times I)$. In this setting, {\it surgery of $F$ along D} refers to the process of removing $\phi(\alpha \times I)$ from $F$, and replacing it with $\phi(D \times \partial I)$.

There are new types of behavior we can describe for $h^{-1}(t)$, as $t$ changes. Much like in the previous section, where we saw 2-sphere components being capped off (or the appearance of 2-sphere components), we may now see disks homotoped to a point on $\partial M$ (or the appearance of new disks). Similarly, we may see a puncture appear or disappear in a component of $h^{-1}(t)$, when $h^{-1}(t)$ moves past a tangency with $\partial M$.

For the remainder of this paper, we will regard the appearance and disappearance of punctures as non-generic, in the following sense: Consider the double of $M$, $DM=M \cup _{\partial M} \overline M$ (where $\overline M$ denotes $M$ with opposite orientation). The function, $h$, also doubles to a function, $Dh:DM \rightarrow I$. At the point where we would see the disappearance of a puncture in a component of $h^{-1}(t)$, we see a compression happen for a component of $Dh^{-1}(t)$. Furthermore, this compression happens at exactly a point of $\partial M$. An arbitrarily small perturbation of $Dh$ makes the compression happen at an interior point of $M$. We now restrict this perturbed version of $Dh$ to $M$, and call the result $h$ again. Where we saw the disappearance of a puncture in a component of $h^{-1}(t)$, we now see a compression, followed by the disappearance of a disk component.

A more significant change in the topology of leaves may now be by $\partial$-compression or $\partial$-decompression (the opposite of a $\partial$-compression). To account for this, we must alter our definition of the complexity of a leaf. Suppose $\mathcal F_t^i$ is a component of $\mathcal F_t$. If $\mathcal F_t^i$ is closed, then define $c(\mathcal F_t^i)$ as before. If $\mathcal F_t^i$ is not closed, then define $c(\mathcal F_t^i)$ to be $0$ if $\mathcal F_t^i$ is a disk, and $1/2-\chi (\mathcal F_t^i)$ otherwise. Let $c(\mathcal F_t)=\sum _i c(\mathcal F_t^i)$, where the sum is taken over all components of $\mathcal F_t$. This measure of complexity will decrease if we see any compression {\it or} $\partial$-compression of $\mathcal F_t$, and it will be $0$ if and only if $\mathcal F_t$ is a collection of spheres and disks.

Let $\{ s_i \}$ be some collection of points in [0,1], such that there is exactly one element of this set between any two consecutive critical values of $h$. Note that we can obtain $\mathcal F_{s_i}$ from $\mathcal F_{s_{i-1}}$ by either adding or removing a 2-sphere or disk, by compressing or de-compressing, or by $\partial$-compressing or $\partial$-decompressing.

Now, let $\{ t_i \}$ be some subcollection of $\{ s_i \}$ such that $\mathcal F_{t_i}$ differs from $\mathcal F_{t_{i+1}}$ by exactly one compression, $\partial$-compression, de-compression, or $\partial$- decompression (and possibly several 2-sphere components and disks). We now define local maxima and minima of $\mathcal F$ precisely as before. That is, a {\it local maximum occurs at} $t_i$ if $c(\mathcal F_{t_i})>c(\mathcal F_{t_{i-1}})$ and $c(\mathcal F_{t_i})>c(\mathcal F_{t_{i+1}})$. As before, if a local maximum (minimum) occurs at $t_i$, then we refer to $\mathcal F_{t_i}$ as a {\it maximal (minimal) leaf} of $\mathcal F$.

Once again we ask: What do the submanifolds of $M$ between consecutive maximal and minimal leaves look like? They are no longer compression bodies. However, by generalizing the definition of a compression body in the appropriate way, we can still give a complete description.

We say a separating surface, $F$, is {\it completely compressible and $\partial$-compressible to one side} if there exists a collection of disjoint compressing disks and $\partial$-compressing disks for $F$ on one side, such that surgery along every disk in this collection yields a collection of spheres which bound balls, or yields a surface which is parallel to some subsurface of $\partial M$.

A $\partial${\it -compression body} is a 3-manifold $W$, equipped with 3 subsurfaces of $\partial W$, which are denoted $\partial _+ W, \partial _- W$, and $\partial _0 W$, such that $\partial _+W$ is completely compressible and $\partial$-compressible, and when compressed, is parallel to $\partial _-W$ (if $\partial _- W \ne \emptyset$), or is a boundary parallel disk, and such that $\partial W=\partial _+ W \cup \partial _- W \cup \partial _0 W$.

        \begin{figure}[htbp]
    \psfrag{0}{$\partial _0 W$}
    \psfrag{-}{$\partial _- W$}
    \psfrag{+}{$\partial _+ W$}
        \vspace{0 in}
        \begin{center}
        \epsfxsize=4 in
        \epsfbox{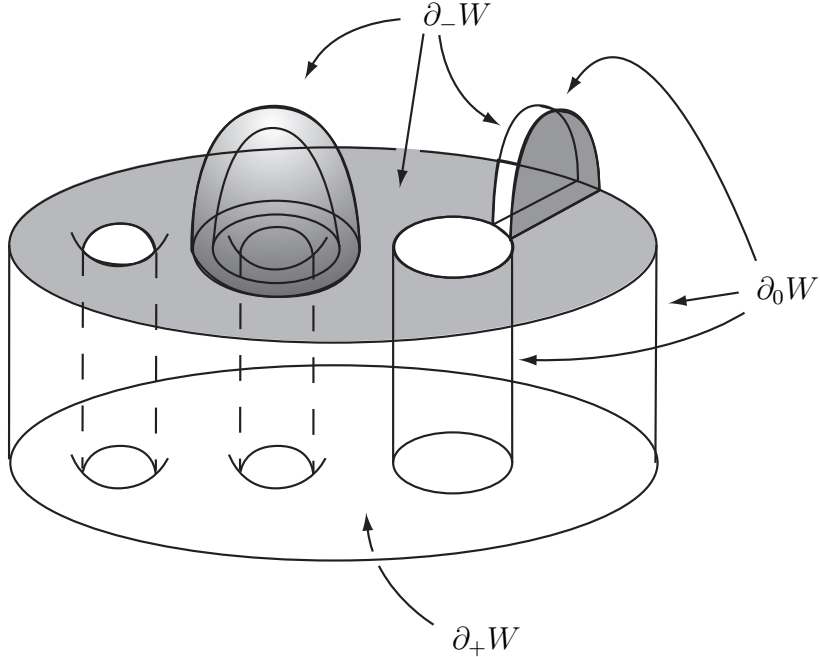}
        \caption{A $\partial$-compression body, made from thickening a twice punctured surface of genus 2 and attaching a 2-handle and a half 2-handle. The shaded regions are $\partial _- W$.}
        \label{f:compbod}
        \end{center}
        \end{figure}

We can also give a constructive description of a $\partial$-compression body, $W$ (see figure \ref{f:compbod}). Let $F$ be some surface, and begin with $F \times I$. Denote $F \times \{0 \}$ by $\partial _+ W$, $(\partial F) \times I$ by $\partial _0 W$, and $F \times \{1\}$ by $\partial _- W$. We now attach a non-empty collection of 2-handles and half 2-handles to $F \times \{1\}$. A {\it half 2-handle} is defined to be $D^2 \times I$, where we think of $\partial D^2=\alpha \cup \beta$, where $\alpha \times I$ is the region we attach to the neighborhood of an arc, $\delta$, in $F \times \{1\}$, such that $\partial \delta \subset \partial (F \times \{1\})$. For each such half 2-handle, we add $D^2 \times \partial I$ to $\partial _-W$, and $\beta \times I$ to $\partial _0 W$. As usual, a 2-handle is just $D^2 \times I$, attached along $(\partial D^2) \times I$. For each such 2-handle added, we add $D^2 \times \partial I$ to $\partial _-W$. Finally, we cap off any 2-sphere components of $\partial _-W$ by 3-balls, and we add any disk components to $\partial _0 W$. Note that $\partial _0 W$ is not in general a product in this setting.

We say that a surface, $F$, in $W$, is $\partial _0${\it -compressible} if there exists a disk, $D$, such that $\partial D=\alpha \cup \beta$, where $D \cap F=\alpha$ is an essential arc on $F$, and $D \cap \partial _0 W=\beta$. One can show that $\partial _-W$ is both incompressible and $\partial _0$-incompressible in $W$.

Another important fact is that $\partial _0 W$ must be incompressible in $W$. To see this, just double $W$ along $\partial _0 W$. Every half 2-handle becomes a 2-handle, so this new manifold is a compression body. A compressing disk for $\partial _0 W$ then doubles to become an essential 2-sphere in a compression body, which cannot happen. We will use this fact in the proof of Theorem \ref{t:bcg}, in Appendix \ref{a:bcg}.

A $\partial${\it -Heegaard splitting} of a manifold, $M$, is a decomposition into two $\partial$-compression bodies, $W$ and $W'$, such that $W \cap W'=\partial _+ W=\partial _+ W'=F$. As before, we denote such a splitting as $W \cup _F W'$. A {\it strongly irreducible $\partial$-Heegaard splitting} is one which has the property that every compressing and $\partial _0$-compressing disk for $F$ in $W$ must have non-empty intersection with every compressing and $\partial _0$-compressing disk for $F$ in $W'$. A {\it quasi-strongly irreducible $\partial$-Heegaard splitting} is one in which every compressing and honest $\partial _0$-compressing disk in $W$ meets every compressing and honest $\partial _0$-compressing disk in $W'$.

We now present an analogous statement to Theorem \ref{t:cg}. First, if $W \cup _F W'$ is a $\partial$-Heegaard splitting of $M$, then let $\partial _-M=\partial _-W \cup \partial _-W'$, and $\partial _0 M=\partial _0 W \cup \partial _0 W'$. A {\it $\partial _0$-compression for $\partial _-M$} is a disk, $D$, such that $\partial D=\alpha \cup \beta$, where $D \cap \partial _-M=\alpha$, $\alpha$ is an essential arc on $\partial _-M$, and $D \cap \partial _0 M=\beta$.

\begin{thm}
\label{t:bcg}
If $W \cup _F W'$ is a quasi-strongly irreducible $\partial$-Heegaard splitting of $M$, then $\partial _-M$ is both incompressible and $\partial _0$-incompressible in $M$.
\end{thm}

We leave the proof of this theorem to Appendix \ref{a:bcg}.

Our picture of a 3-manifold with boundary is now completely analogous to the last section. A Morse function on $M$ which induces a Morse function on $\partial M$ defines a singular foliation, $\mathcal F$. The minimal leaves of $\mathcal F$ break up $M$ into submanifolds, each one having a $\partial$-Heegaard splitting surface which is a maximal leaf of $\mathcal F$.

\section{Examples}
\label{ex}

\begin{example}
A quasi-strongly irreducible $\partial$-Heegaard splitting may not be strongly irreducible.
\end{example}

\begin{figure}[htbp]
    \psfrag{W}{$W$}
    \psfrag{w}{$W'$}
    \psfrag{z}{$\partial _0 W$}
    \psfrag{0}{$\partial _0 W'$}
    \psfrag{-}{$\partial _- W$}
        \vspace{0 in}
        \begin{center}
        \epsfxsize=4 in
        \epsfbox{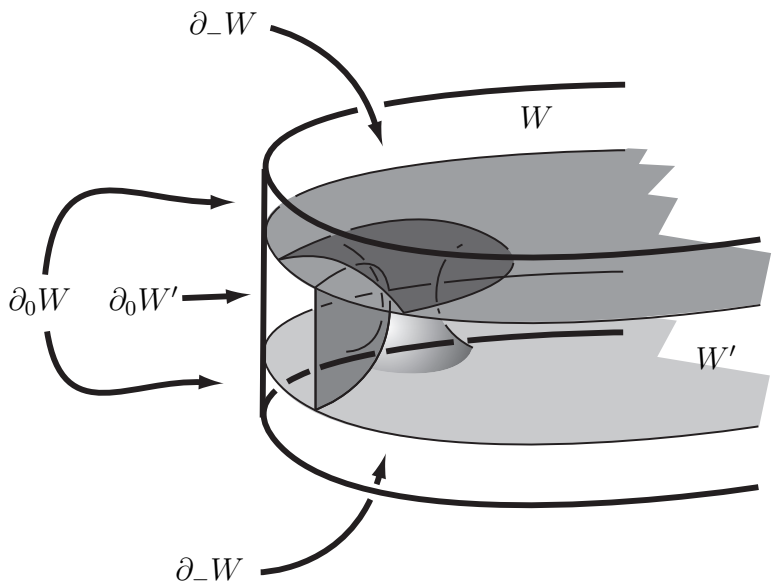}
        \caption{Disjoint $\partial$-compressions on opposite sides of $F$, where one is not honest.}
        \label{f:sigma}
        \end{center}
        \end{figure}

\begin{proof}
Let $M=\Sigma \times I$, where $\Sigma$ is some surface with nonempty boundary, other than $D^2$. Let $F$ be the surface obtained by connecting $\Sigma \times \{1/3\}$ to $\Sigma \times \{2/3\}$ by an unknotted, boundary compressible tube, as in figure \ref{f:sigma}. Let $W$ be the side of $F$ which contains $\Sigma \times \{0\}$ and $\Sigma \times \{1\}$, and let $W'$ denote the other side of $F$. Then $W \cup _F W'$ is a $\partial$-Heegaard splitting, where $\partial _-W = \Sigma \times \{0\} \cup \Sigma \times \{1\}$, $\partial _0 W = \partial \Sigma \times \{[0,1/3] \cup [2/3,1]\}$, $\partial _0 W'=\partial \Sigma \times [1/3,2/3]$, and $\partial _- W' =\emptyset$. $F$ is a quasi-strongly irreducible $\partial$-Heegaard splitting of $M$, but is not strongly irreducible.
\end{proof}

\begin{example}
\label{e:bridgeknot}
Knots and links in bridge position yield $\partial$-Heegaard splittings.
\end{example}

\begin{proof}
Consider a knot (or link), $K \subset S^3$, which is in bridge position. That is, there is some height function, $h$, on $S^3$, in which all of the minima of $K$ are below all of the maxima. Suppose $S=h^{-1}(1/2)$ is a level 2-sphere which separates the minima from the maxima. Let $M^K$ denote $S^3$ with a neighborhood of $K$ removed, and $S^K=S \cap M^K$. If $W$ is the region of $M^K$ above $S^K$, then $W$ is a $\partial$-compression body, where $\partial _+ W =S^K$, $\partial _-W=\emptyset$, and $\partial _0 W$ is the remainder of $\partial W$. Likewise, the region of $M^K$ below $S^K$ is a $\partial$-compression body, and so $S^K$ is a $\partial$-Heegaard surface.
\end{proof}

Suppose that $K \subset S^3$ is an arbitrary knot or link, and $h$ is some height function on $S^3$, which is a Morse function when restricted to $K$. Let $\{ q_j \}$ denote the critical values of $h$ restricted to $K$, and let $q'_j$ be some point in the interval $(q_j, q_{j+1})$. Then the {\it width} of $K$ is the sum over all $j$ of $|h^{-1}(q'_j) \cap K|$. If $K$ realizes its minimal width, then we say $K$ is in {\it thin position} (see \cite{gabai:87}).

\begin{example}
\label{e:thinknot}
Knots and links in thin and bridge position yield strongly irreducible $\partial$-Heegaard splittings.
\end{example}

\begin{figure}[htbp]
    \psfrag{f}{$\delta$}
    \psfrag{d}{$\delta '$}
    \psfrag{e}{$\delta ''$}
    \psfrag{A}{$A$}
    \psfrag{G}{$D$}
    \psfrag{S}{$A'$}
    \psfrag{Q}{$D'$}
    \psfrag{Z}{$\partial _0 D$}
    \psfrag{W}{$\partial _- D$}
    \psfrag{B}{$\partial _0 M$}
    \psfrag{F}{$F \times \{1\}$}
    \psfrag{a}{$\alpha$}
    \psfrag{b}{$\beta$}
    \psfrag{V}{$W$}
    \psfrag{X}{$S^K$}
        \vspace{0 in}
        \begin{center}
        \epsfxsize=3.5 in
        \epsfbox{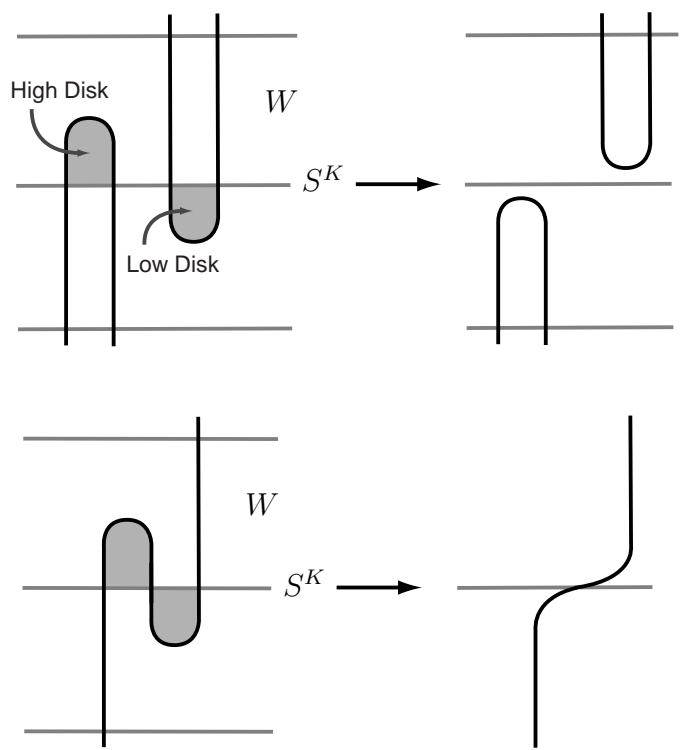}
        \caption{Isotoping $K$ when there are disjoint High and Low Disks.}
        \label{f:High&Low}
        \end{center}
        \end{figure}

\begin{proof}
Suppose that the knot (or link), $K$, of Example \ref{e:bridgeknot} is in thin position, as well as bridge position. Also, assume $K$ is not the unknot. We will depart from standard terminology a bit here. A $\partial$-compressing disk for $S^K$ which lies entirely above it will be referred to as a {\it ``high disk"}, and one which lies below it will be called a {\it ``low disk"}. If we see a high disk which is disjoint from a low disk, then we can isotope $K$ as in figure \ref{f:High&Low}, to obtain a presentation of smaller width. Hence, any $\partial$-compression above $S^K$ must intersect every $\partial$-compression below it. Now, suppose there is a compressing disk, $D$, for $S^K$ in $W$. Then $D$ caps off some maxima of $K$, all of which correspond to high disks. Also, since $D$ is a compressing disk for $S^K$, there must be some maxima of $K$ (and hence, some high disks) on the other side of $D$ in $W$. Similarly, any compressing disk, $D'$, for $S^K$ which lies below it must have low disks on both sides. If $D \cap D' =\emptyset$, then we can conclude that there were disjoint high and low disks, and hence, $K$ was not thin. Likewise, it is easy to rule out the case where we have a compressing disk for $S^K$ on one side, which is disjoint from a $\partial$-compression on the other side. Our conclusion is that $S^K$ is strongly irreducible.
\end{proof}

\begin{dfn}
A $\partial$-Heegaard splitting, $W \cup _F W'$ is {\it stabilized} is there exist compressing disks on each side of $F$ which meet in a single point, {\it or} a $\partial$-compressing disk on one side that meets a compressing disk on the other in a single point.
\end{dfn}

\bigskip

{\bf Exercise.}
{\it A stabilized $\partial$-Heegaard splitting, $W \cup _F W'$, either fails to be strongly irreducible, is the genus 1 Heegaard splitting of $S^3$, or $F$ is an unknotted annulus in $B^3$.}

\medskip

It is interesting to note that the double of an unknotted annulus in $B^3$ gives the genus 1 splitting of $S^3$. For the remainder of this paper, we shall always assume that all strongly irreducible $\partial$-Heegaard splittings are not stabilized, whereas quasi-strongly irreducible splittings may be stabilized. In light of the above exercise, this does not greatly reduce possible applications.

\section{Mini-Lmax Foliations}
\label{s:closed}

Let $M$ be a compact, orientable, irreducible 3-manifold, and suppose $h:M \rightarrow [0,1]$ is a Morse function, where we do not require that $\partial M \subset h^{-1}(0) \cup h^{-1}(1)$. Let $\mathcal F$ be the singular foliation induced by $h$, and let $Lmax(\mathcal F)=\{ c(\mathcal F_{t_i})$ such that a local maximum occurs at $t_i \}$ (where we include repeated integers). We arrange this set in non-increasing order, and compare two such sets lexicographically. This gives us a way of comparing two singular foliations of $M$.

\begin{dfn}
$\mathcal F$ is a {\it mini-Lmax foliation} if for every foliation, $\mathcal F'$, of $M$, $Lmax(\mathcal F)\le Lmax(\mathcal F')$.
\end{dfn}

The reason for the terminology is that this is a strict generalization of the concept of $\mathcal F$ being minimax (see, for example, \cite{rubinstein:93}). The number we first want to minimize under this definition is the maximal value of $c(\mathcal F_t)$. Hence, if $\mathcal F$ is mini-Lmax, then $\mathcal F$ is minimax. Now, among all such foliations, choose the subset such that the second largest value of $c(\mathcal F_t)$ is minimal. If we repeat this process, we arrive at the set of mini-Lmax foliations.

This definition is also extremely similar to the complexity defined in \cite{st:94}, the only difference being that in that paper, the sets which one compares consist of all values of $c(\mathcal F_{t_i})$, rather than just the maximal values, and the requirement is made that $\partial M \subset h^{-1}(0) \cup h^{-1}(1)$.

\begin{thm}
\label{bmaxima}
Let $\mathcal F$ be a mini-Lmax foliation of $M$. Then the maximal leaves of $\mathcal F$ are strongly irreducible $\partial$-Heegaard surfaces for the submanifolds obtained by cutting $M$ along minimal leaves.
\end{thm}

\begin{proof}
Recall from \cite{st:94} that the analogous theorem was true because if we ever saw a compression on the ``top" side of a maximal leaf, that was disjoint from a compression on the ``bottom" side, then we could decompress along the upper one before compressing along the lower one. This gives rise to a foliation of the same manifold with lower $Lmax(\mathcal F)$.

The situation is precisely the same here. If we see a boundary compression on one side which is disjoint from either a compression or another boundary compression on the other, then we can re-arrange the order of compressions, de-compressions, $\partial$-compressions and $\partial$-decompressions to obtain a foliation with smaller $Lmax(\mathcal F)$.
\end{proof}

If $\mathcal F$ is any foliation which satisfies the conclusion of Theorem \ref{bmaxima}, then we say $\mathcal F$ is {\it locally mini-Lmax}. In fact, we shall even refer to $\mathcal F$ as locally mini-Lmax if the maximal leaves are only quasi-strongly irreducible Heegaard surfaces. Note that a strongly irreducible ($\partial$-)Heegaard splitting of any manifold gives rise to an example of a locally mini-Lmax foliation, since any Heegaard surface can be realized as the maximal leaf in a singular foliation with only one maximal leaf. (In fact, we can take this as the definition of a Heegaard surface).

Theorem \ref{bmaxima} gives a very nice picture of a manifold with boundary. In particular, we see that any manifold that admits a locally mini-Lmax foliation can be decomposed into two sets of $\partial$-compression bodies, $\{W_i\}$, and $\{W'_i\}$, where $\partial _+ W_i =\partial _+ W'_i$, and $\partial _- W'_i = \partial _- W_{i+1}$. Also, if $1 \le i \le n$, then $\partial M=\partial _- W_1 \cup (\cup \partial _0 W_i) \cup (\cup \partial _0 W'_i) \cup \partial _- W'_n$. Let $\partial _0 M=(\cup \partial _0 W_i) \cup (\cup \partial _0 W'_i)$. We now immediately deduce the following theorem.

\begin{thm}
\label{bminima}
If $\mathcal F$ is a locally mini-Lmax foliation of $M$, then the minimal leaves of $\mathcal F$ are incompressible and $\partial _0$-incompressible in $M$.
\end{thm}

\begin{proof}
Theorem \ref{t:bcg} implies that the minimal leaves of $\mathcal F$ are incompressible and $\partial _0$-incompressible in the submanifolds obtained by cutting $M$ along minimal leaves. A standard innermost disk/outermost arc argument shows they are incompressible and $\partial _0$-incompressible in $M$.
\end{proof}

As an immediate corollary, we obtain:

\begin{cor}
Let $\mathcal F$ be a locally mini-Lmax foliation of $M$. If $\partial M=\partial _0 M$, or if $\partial _0 M \cap \partial _-M=\emptyset$, then the minimal leaves of $\mathcal F$ are incompressible and $\partial$-incompressible in $M$.
\end{cor}

\section{Foliations and Knots and Links}
\label{s:knots}

We would now like to discuss further singular foliations in the complement of knots and links. Suppose $(K,\partial K) \subset (M, \partial M)$ is an embedded 1-manifold. Let $M^K$ denote $M$, with a small neighborhood of $K$ removed. If $X$ is some subset of $M$, then let $X^K=X \cap M^K$.

\begin{dfn}
\label{d:locunknot}
A 1-manifold $(K,\partial K) \subset (M,\partial M)$ is {\it locally tangled} if there is a ball, $B \subset M$, such that $(\partial B)^K$ is incompressible and $\partial$-incompressible in $M^K$, or such that $K \subset B$. If no such ball exists, then $K$ is {\it locally untangled}.
\end{dfn}

For the remainder of this paper, we will assume that $M^K$ is irreducible. If $\mathcal F$ is a singular foliation of $M$ arising from some height function, then let $\mathcal F^K=\mathcal F \cap M^K$. A leaf of $\mathcal F^K$ shall be denoted as $\mathcal F_t^K$.

\begin{dfn}
\label{d:miniLmaxK}
Suppose $\mathcal F$ is a singular foliation of $M$. A 1-manifold, $K$, is in a position which is {\it mini-Lmax with respect to $\mathcal F$} (or simply {\it mini-Lmax}, when it is clear what $\mathcal F$ is), if $\partial K \subset h^{-1}(0) \cup h^{-1}(1)$, and $K$ cannot be isotoped to reduce $Lmax(\mathcal F^K)$.
\end{dfn}

We are now in a position to generalize Example \ref{e:thinknot}.

\begin{thm}
\label{maxwrtk}
Let $M$ be an irreducible 3-manifold other than $S^3$. Suppose $K$ is a locally untangled 1-manifold, which is mini-Lmax with respect to a locally mini-Lmax foliation, $\mathcal F$, such that no component of $K$ can be isotoped onto a leaf of $\mathcal F$. Then the maximal leaves of $\mathcal F^K$ are quasi-strongly irreducible $\partial$-Heegaard surfaces for the submanifolds of $M^K$ that arise when we cut along the minimal leaves.
\end{thm}

\begin{proof}
Suppose $P$, $Q$, and $Q_*$ are leaves of $\mathcal F$ such that $P^K$ is a maximal leaf of $\mathcal F^K$, and $Q_*^K$ and $Q^K$ are consecutive minimal leaves of $\mathcal F^K$ which ``sandwich" $P^K$ (one or both may be empty). Let $W$ be the region of $M$ between $P$ and $Q$ (see figure \ref{f:W-K}), and $W_*$ be the region between $P$ and $Q_*$.

\medskip

{\bf Case 1.} There are compressing (or honest $\partial$-compressing) disks, $D$ and $D_*$, for $P^K$, in $W^K$ and $W_*^K$, which are compressing (or honest $\partial$-compressing) disks for $P$. Then not only is $P^K$ a maximal leaf for $\mathcal F^K$, but also $P$ is a maximal leaf of $\mathcal F$. If $\partial D \cap \partial D_* =\emptyset$, then $P$ fails to be a strongly irreducible Heegaard surface, and hence, $\mathcal F$ is not locally mini-Lmax. Since the local mini-Lmaximality of $\mathcal F$ was a hypothesis of the Theorem, this is a contradiction.

        \begin{figure}[htbp]
        \vspace{0 in}
        \begin{center}
        \epsfxsize=2.5 in
        \epsfbox{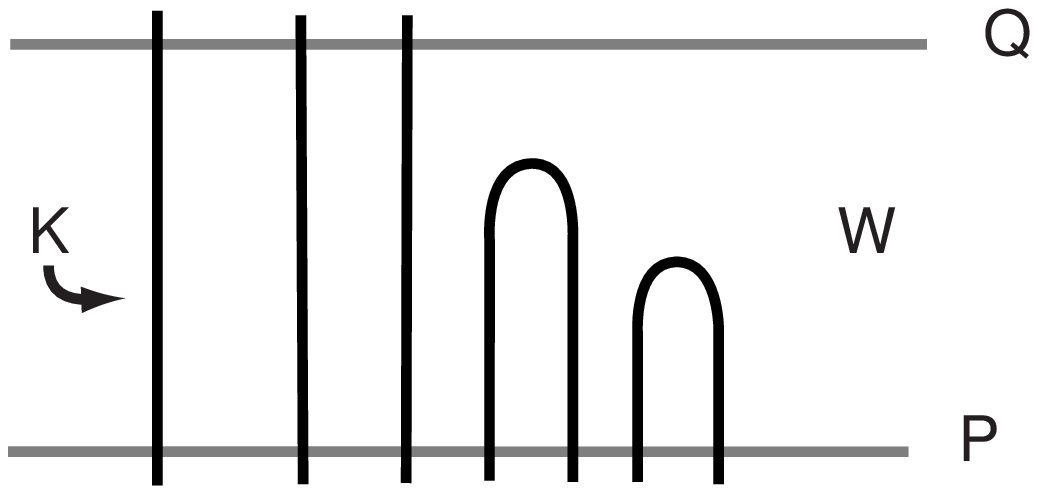}
        \caption{$W^K$.}
        \label{f:W-K}
        \end{center}
        \end{figure}

{\bf Case 2.} Suppose $D$ and $D_*$ are disjoint compressing (or honest $\partial _0$-compressing) disks for $P^K$, but not $P$. Then we are in a very similar situation to Example \ref{e:thinknot}. If they are both honest $\partial _0$-compressing disks, then since no component of $K$ can be isotoped onto a leaf of $\mathcal F$, we can do one of the moves depicted in Figure \ref{f:High&Low} to reduce $Lmax(\mathcal F^K)$ (if some component of $K$ could be isotoped onto a leaf, we'd have another possibility to consider, whose effect would be such an isotopy). Note that this is the only place in the proof of this Theorem where we use the assumption of honesty. This is necessary because an honest $\partial _0$-compressing disk must look like a high or low disk, whereas there may be many possibilities for a $\partial$-compressing disk which is not honest.

If both $D$ and $D_*$ are compressing disks, then $\partial D$ bounds a disk, $E$, on $P$, which must be punctured by $K$. Since $Q ^ K$ is the first minimal leaf after $P^K$, $K \cap W$ must consist of a collection of vertical arcs and trivial arcs, which contain a single maximum, as in Figure \ref{f:W-K}. Since $D$ lies in $W ^ K$, we must see an arc of the later type in the ball bounded by $D \cup E$. Such arcs always co-bound high disks. Similarly, $\partial D_*$ bounds a disk, $E_*$, on $P$, and we see low disks inside the ball bounded by $D_* \cup E_*$. If $\partial D \cap \partial D_*=\emptyset$, then there are two cases.

\smallskip

{\bf Subcase 2.1.}
If $E \cap E_*=\emptyset$, then we see disjoint high and low disks for $P$, which is again a contradiction.

\smallskip

        \begin{figure}[htbp]
        \vspace{0 in}
        \begin{center}
        \epsfxsize=4.5 in
        \epsfbox{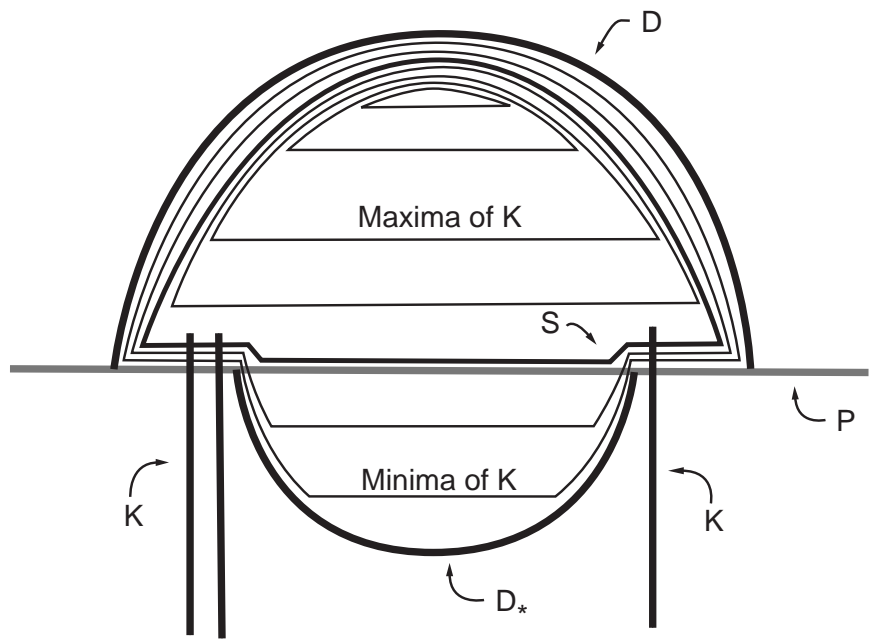}
        \caption{The foliation, $\mathcal F^B$, of $B$.}
        \label{f:Foliation}
        \end{center}
        \end{figure}

{\bf Subcase 2.2.} The other case is when $E_* \subset E$ (or $E \subset E_*$; the proof will be symmetric). Let $B$ denote the union of the ball bounded by $D \cup E$, and the ball bounded by $D_* \cup E_*$. We now claim that $(\partial B)^K$ is incompressible and $\partial$-incompressible in $B^K$.

Consider the foliation, $\mathcal F^B$, of $B$ depicted in Figure \ref{f:Foliation}. For each leaf, $\mathcal F_t$, of $\mathcal F$, which intersects $\partial B$, we construct a leaf, $\mathcal F_t^B$, of $\mathcal F^B$ as follows: let $(\partial B)^+$ denote the subset of $\partial B$ above $\mathcal F_t$. Now, let $\mathcal F_t^B = (\mathcal F_t \cap B) \cup (\partial B)^+$. This leaf can be pushed slightly into $B$, so that the foliation, $\mathcal F^B$, is well defined over all of the interior of $B$. To complete $\mathcal F^B$, we simply add a leaf which is precisely $\partial B$. Away from a neighborhood of the boundary of $B$, this foliation just looks like $\mathcal F$.

Note that $K \cap B$ is in bridge position with respect to $\mathcal F^B$. We now focus on the 2-sphere, $S$, depicted in figure \ref{f:Foliation} (which is precisely the leaf $\mathcal F_t^B$, of $F^B$, where $\mathcal F_t=P$). Note that every honest $\partial _0$-compressing disk for $S^K$ in $B^K$ which is on one side of $S^K$ (i.e. every high disk) must intersect every honest $\partial _0$-compressing disk for $S^K$ on the other side (i.e. every low disk). If not, then we could isotope $K$ inside $B$, and decrease $Lmax(\mathcal F^K)$, a contradiction. We are now in precisely the same situation as in Example \ref{e:thinknot}, so we may conclude that $S$ is a quasi-strongly irreducible $\partial$-Heegaard surface for $B^K$. Hence, by Theorem \ref{t:bcg}, $(\partial B)^K$ is incompressible and $\partial _0$-incompressible in $B^K$. We now compress and $\partial$-compress $(\partial B)^K$ completely to the outside of $B$ to obtain a sphere, $S'$, which bounds a ball in $M$ (by irreducibility), which either contains $K$, or such that $(S')^K$ is incompressible and $\partial$-incompressible in $M^K$. (Note that $\partial B$ cannot compress away to nothing outside $B$, since $M$ is not homeomorphic to $S^3$). This shows that $K$ was locally tangled, violating the hypothesis of Theorem \ref{maxwrtk}.

\smallskip

Similarly, if $D$ is an honest $\partial _0$-compression and $D_*$ is a compression, then we can find disjoint high and low disks for $P$, or show $K$ was locally tangled.

\medskip

{\bf Case 3.} The last case we need to consider is when $D$ and $D_*$ are compressing (or $\partial _0$-compressing) disks for $P^K$, but only $D_*$ is a compressing (or $\partial _0$-compressing) disk for $P$. In this case, as in the preceding case, we see a high disk, $H \subset W$, such that $\partial H=\alpha \cup \beta$, where $H \cap K =\beta$, $H \cap P=\alpha$, and $\partial D_* \cap \alpha = \emptyset$. This situation, too, never occurs for a maximal leaf in a mini-Lmax foliation. We simply compare this foliation to the one isotopic to $\mathcal F$, where we pass through the maxima of $K \cap H$ {\it before} decompressing along $D_*$. In other words, we can reduce $Lmax(\mathcal F^K)$ by using $H$ to isotope $K$ below $P$.

\smallskip

In short, we have shown that if $D$ is any compressing (or honest $\partial _0$-compressing) disk for $P^K$ in $W^K$, and $D_*$ is a compressing (or honest $\partial _0$-compressing) disk for $P^K$ in $W_*^K$, then $\partial D \cap \partial D_* \ne \emptyset$. Hence, $P^K$ is a quasi-strongly irreducible $\partial$-Heegaard surface for $(W^K) \cup (W_*^K)$.
\end{proof}

\begin{thm}
\label{minwrtk}
For $M, K$, and $\mathcal F$ as in the statement of Theorem \ref{maxwrtk}, the minimal leaves of $\mathcal F^K$ are incompressible and $\partial _0$-incompressible in $M^K$.
\end{thm}

\begin{proof}
As in Theorem \ref{bminima}, an application of Theorem \ref{t:bcg} tells us that $Q^K$ and $Q_*^K$ are incompressible and $\partial _0$-incompressible in $(W^K) \cup (W_*^K)$, and a standard innermost disk/ outermost arc argument shows they are incompressible and $\partial _0$-incompressible in $M^K$.
\end{proof}

For any triple, $(M, K, \mathcal F)$, which satisfies the conclusion of Theorem \ref{maxwrtk}, we say $K$ is {\it locally mini-Lmax with respect to $\mathcal F$}, or, when it is clear, just {\it locally mini-Lmax}. Note that the local mini-Lmaximality of $K$ is sufficient to prove Theorem \ref{minwrtk}.

We can make this condition a bit easier to state if we alter our language a bit. {\it For the remainder of this paper, we shall refer to ANY compressing or $\partial _0$-compressing disk for $P^K$ in $W^K$ as a high disk, and ANY compressing or $\partial _0$-compressing disk for $P^K$ in $W_*^K$ as a low disk.} Now, the condition that $K$ is locally mini-Lmax with respect to $\mathcal F$ means that we have no disjoint high and low disks.

We conclude this section with a generalization of the main result of \cite{thompson:97}. Suppose $K$ is some knot embedded in $M$, and $H$ is a strongly irreducible Heegaard surface in $M$. Let $h:M \rightarrow [0,1]$ be a Morse function such that $H$ is the maximal leaf of the singular foliation induced by $h$. If the maxima of $K$ are all above $H$, and the minima all below, then we say $K$ {\it is in bridge position with respect to} $H$.

\begin{thm}
\label{t:generalthompson}
Let $H$ be a strongly irreducible Heegaard splitting of a closed, orientable, irreducible 3-manifold, $M$ (If $M$ is homeomorphic to $S^3$, then let $H$ be any embedded 2-sphere). Let $h:M \rightarrow [0,1]$ be a Morse function such that $H$ is the maximal leaf of the singular foliation, $\mathcal F$, induced by $h$. Let $K$ be any 1-manifold embedded in $M$, which has no component isotopic onto $H$, such that $M^K$ is irreducible. If $K$ is mini-Lmax with respect to $\mathcal F$, then either $K$ is in bridge position with respect to $H$, or there is a meridional, incompressible, $\partial$-incompressible surface in $M^K$, which has genus less than or equal to that of $H$.
\end{thm}

\begin{proof}
If $M$ is homeomorphic to $S^3$, then this is precisely the main result of \cite{thompson:97}. So, assume $M$ is not $S^3$, and $K$ is mini-Lmax with respect to $\mathcal F$. First, if $K$ is locally tangled, then by definition the Theorem is true. So, assume $K$ is not locally tangled. If $K$ is not in bridge position with respect to $H$, then there is some minimal leaf of $\mathcal F^K$. Theorem \ref{minwrtk} now implies that this surface is incompressible and $\partial$-incompressible in $M^K$. Since $H$ is the maximal leaf of $\mathcal F$, every leaf of $\mathcal F^K$ has genus less than or equal to that of $H$.
\end{proof}

\section{Normal Surfaces: Definitions}
\label{s:normal}

\begin{figure}[htbp]
        \vspace{0 in}
        \begin{center}
        \epsfxsize=3.25 in
        \epsfbox{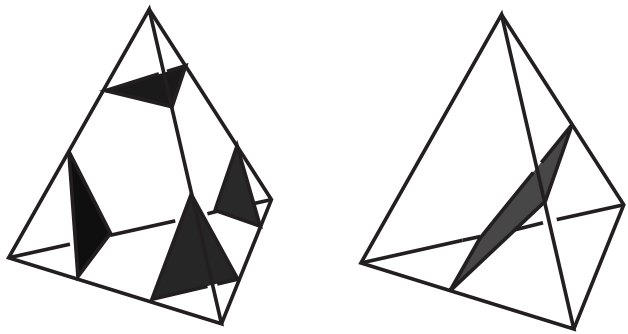}
        \caption{Normal Disks.}
        \label{f:Normal}
        \end{center}
        \end{figure}

In this section, we discuss the necessary background material on normal surfaces. A {\it normal curve} on the boundary of a tetrahedron is a simple loop which is transverse to the 1-skeleton, made up of arcs which connect distinct edges of the 1-skeleton. The {\it length} of such a curve is simply the number of times it crosses the 1-skeleton. A {\it normal disk} in a tetrahedron is any embedded disk, whose boundary is a normal curve of length three or four, and whose interior is contained in the interior of the tetrahedron, as in figure \ref{f:Normal}.

A {\it normal surface} in $M$ is the image of an embedding, $p$, of some surface, $(F,\partial F)$, into $(M,\partial M)$, such that $p(F)$ is a union of normal disks. In addition, we say $p(F)$ is an {\it almost normal surface} if it consists of all normal disks, plus one additional piece in one tetrahedron. This piece can be either a disk with normal boundary of length 8 (depicted in figure \ref{f:Octagon}), two normal disks connected by a single unknotted tube (as in figure \ref{f:1Tube}), or two normal disks connected by a band along $\partial M$ (see figure \ref{f:halftube}). The first two types of almost normal surfaces were first explored by Rubinstein in \cite{rubinstein:93}, and later used by Thompson \cite{thompson:94} and Stocking \cite{stocking:96}. This paper generalizes many of those results to surfaces of the third type.

\section{Normal and Almost Normal Surfaces and Mini-Lmax Foliations}
\label{s:nlan}

One application of the results we have discussed thus far comes about when we let $K$ be the 1-skeleton of a pseudo-triangulation of $M$. To make this more precise, suppose $T$ is any pseudo-triangulation of $M$ (i.e. an expression of $M$ as a union of 3-simplices, where any two such 3-simplices intersect in a (possibly empty) collection of lower dimensional simplices). Let $T_n$ denote the $n$-skeleton of $T$. We now focus on singular foliations which arise from height functions, as before. However, we must make a few  additional restrictions: for $M$ closed, we require that $T_0$ consists of a single vertex. If $\partial M \ne \emptyset$, then we require that $T_0 \subset \partial M$, and that each component of $\partial M$ contains exactly one component of $T_0$. In either case, we also need that $T_0 \subset h^{-1}(0) \cup h^{-1}(1)$. In addition, we require that the only normal 2-sphere in $M$ (if any) is a link of $T_0$. Finding such a triangulation is essentially the first step in the original proofs of the results of this section (\cite{rubinstein:93}, \cite{thompson:94}, \cite{stocking:96}), and we find it necessary in our approach as well. The proof that any irreducible manifold admits such a triangulation can be found in \cite{jr:98}. The reason here for this assumption is that by \cite{jr:88}, we know that any normal 2-sphere is incompressible in the complement of $T_1$. So, if there is a non-vertex linking normal 2-sphere, then $T_1$ is locally tangled, and hence, we will not be able to apply Theorem \ref{maxwrtk}.

Furthermore, in order to make sense of the definitions given in the previous sections, we must push the interiors of the edges of $T_1$ which lie on $\partial M$ slightly into $M$, as well as the interiors of the boundary 2-simplices. The reason for this is that if we see a leaf of $\mathcal F$ become tangent to an edge of $T_1$ which lies on $\partial M$, and then pass through it, we would like to say $c(\mathcal F_t^{T_1})$ has changed. Alternatively, we could have originally defined $c(\mathcal F_t^{T_1})$ to be $c(\mathcal F_t)+|\mathcal F_t \cap T^1|$. Had we done this, all of the results of the preceding sections would have been the same. For closed manifolds, this complexity is {\it exactly} the same as the complexity we originally used.

\begin{dfn}
\label{bubble}
Suppose $\mathcal F_t$ is a leaf of $\mathcal F$ in $M$. A {\it bubble} for $\mathcal F_t$ is a ball, $B$, such that $\partial B= D_1 \cup D_2$, where $D_1$ and $D_2$ are disks, $D_1$ is contained in a single tetrahedron, $\mathcal F_t \cap B=D_2$, $D_2 \cap T^2 \ne \emptyset$, and $D_2 \cap T_1 = \emptyset$.
\end{dfn}

\begin{lem}
\label{nobubbles}
Suppose $T_1$ is mini-Lmax with respect to a locally mini-Lmax foliation, $\mathcal F$. Given some finite collection of non-parallel leaves (i.e. the subset of $\mathcal F^{T_1}$ between any two consecutive leaves of this collection is not a product foliation), we may isotope $\mathcal F$ to obtain a foliation in which no leaf in this collection has any bubbles, and in which $T_1$ is still mini-Lmax with respect to $\mathcal F$.
\end{lem}

\begin{proof}
Suppose $B$ is a bubble for $\mathcal F_t$, where $\partial B=D_1 \cup D_2$, as in Definition \ref{bubble}. We can use $B$ to guide an isotopy from $D_2$ to $D_1$. This may push other leaves which had non-empty intersection with $int(B)$, but it can only destroy bubbles for those leaves, too. Also, the isotopy leaves behind a ``hole" in its wake, but it is easy to fill in intermediate leaves to complete the foliation of $M$. Note that the leaves which we fill in are all parallel to the one just isotoped, so we have not affected any other leaf in our collection. The isotopy is supported on a neighborhood of $B$, which is disjoint from $T_1$. Hence, if $T_1$ was locally minimax with respect to $\mathcal F$, then so is our new foliation. Since there are a finite number of leaves in our collection, and a finite number of bubbles for each, we arrive at a foliation with the desired properties.
\end{proof}

\begin{dfn}
\label{complete}
A {\it complete collection of minimal (maximal) leaves} for $\mathcal F$ is a finite collection, $\{ \mathcal F_{t_i} \}$, such that for every minimal (maximal) leaf, $\mathcal F_t$, of $\mathcal F$, there is an $i$ such that the foliation between $\mathcal F_t$ and $\mathcal F_{t_i}$ is a product. Similarly, a {\it complete collection of minimal (maximal) leaves for $\mathcal F^K$} is a finite collection, $\{ \mathcal F_{t_i} \}$, such that for every minimal (maximal) leaf, $\mathcal F_t ^ K$, of $\mathcal F^K$, there is an $i$ such that the foliation between $\mathcal F_t ^ K$ and $\mathcal F_{t_i} ^ K$ is a product in $M ^ K$.
\end{dfn}

\begin{thm}
\label{normal}
Suppose $\mathcal F$ is a locally mini-Lmax foliation of $M$, and $T_1$ is mini-Lmax with respect to $\mathcal F$. Then we may isotope $\mathcal F$, keeping $T_1$ mini-Lmax, so that every leaf of a complete collection of minimal leaves for $\mathcal F^{T_1}$ is a normal surface.
\end{thm}

\begin{proof}
Let $\{ \mathcal F_{t_i} \}$ be a complete collection of minimal leaves for $\mathcal F^{T_1}$. We begin by using Lemma \ref{nobubbles} to isotope $\mathcal F$ so there are no bubbles for any leaf in this collection.

Now, let $\mathcal F_t$ be a leaf in our collection, let $\tau$ be some tetrahedron in $T$, and let $\Delta$ be a face of $\tau$. First, we examine the possibilities for $\mathcal F_t \cap \Delta$. Let $\gamma$ be an innermost simple closed curve, bounding a disk, $D_1$ in $\Delta$. By Theorem \ref{minwrtk}, $\gamma$ must bound a disk, $D_2$, in $\mathcal F_t^{T_1}$. $M^{T_1}$ is irreducible (it's a handle-body), so $D_1 \cup D_2$ bounds a bubble for $\mathcal F_t$. This is a contradiction, so we see no simple closed curves in any face.

If there are any curves which run from one edge of $\Delta$ to itself, then there is an outermost such one. Let $D$ denote the sub-disk it cuts off in $\Delta$. Then $D$ is a $\partial _0$-compressing disk for $\mathcal F_t^{T_1}$, also contradicting Theorem \ref{minwrtk}. We conclude that $\mathcal F_t \cap \Delta$ is a collection of normal arcs.

We now consider the possibilities for $\mathcal F_t \cap \partial \tau$. It is easy to show that the only possibilities for normal loops are curves of length 3, or 4n (see, for example, \cite{thompson:94}). If there are any curves of length greater than 4, then there must be a disk, $D$, such that $\partial D=\alpha \cup \beta$, where $D \cap T_1 =\alpha$, and $D \cap \mathcal F_t =\beta$ (see \cite{thompson:94}). This is a $\partial$-compressing disk for $\mathcal F_t^{T_1}$, which is again a contradiction. We conclude that $\mathcal F_t \cap \partial \tau$ consists of normal loops of length 3 and 4.

Finally, it follows from Theorem \ref{minwrtk} that every loop of $\mathcal F_t \cap \partial \tau$ bounds a disk on $\mathcal F_t^{T_1}$. Since we have already ruled out simple closed curves in faces of $\tau$, such disks must lie entirely inside $\tau$. We conclude $\mathcal F_t$ is a normal surface.
\end{proof}

Our goal now is to show that once bubbles are removed from maximal leaves, they become almost normal in $M$. First, we shall need a few lemmas.

\begin{lem}
\label{nlarcs}
Suppose $\mathcal F_{t_1} ^ {T_1}$ and $\mathcal F_{t_2} ^ {T_1}$ are consecutive singular leaves of $\mathcal F ^ {T_1}$ such that for each $t \in (t_1,t_2)$, $\mathcal F_t$ is a maximal leaf. If $\partial M=\emptyset$, then there exists a $t' \in (t_1,t_2)$ such that for every 2-simplex, $\Delta$, in $T$, $\mathcal F_{t'} \cap \Delta$ is a collection of normal arcs, and simple closed curves which are inessential on $\mathcal F_{t'} ^ {T_1}$. If $\partial M \ne \emptyset$, then we also allow either a single non-normal arc, or two non-normal arcs on distinct edges of some 2-simplex, which lies on the boundary of $M$.
\end{lem}

\begin{proof}
This proof is similar to many standard arguments which use thin position, but the main idea is taken from \cite{thompson:94}, Claim 4.4. The first step is to consider the nature of the singularities at $t_1$ and $t_2$. If we see a minimum of $T_1$ at $t_1$, and a maximum at $t_2$, then we are in precisely the situation described in \cite{thompson:94}, Claim 4.4. We include the proof here for completeness. Just after $t_1$, we see a low disk for $\mathcal F_t$, contained in the face of some tetrahedron. Similarly, just before $t_2$, we see a high disk in some face. But we never see low disks and high disks at the same time in a maximal leaf, which are disjoint or intersect in a single point. We conclude that there must be some intermediate value where there are no high or low disks in the faces of any tetrahedron, completing the proof in this case. Note that an innermost simple closed curve, which is essential in $\mathcal F_t ^ {T_1}$, bounds a subdisk of $\Delta$ which must be a high or a low disk (recall our modified definition of high and low disks, given just after the proof of Theorem \ref{minwrtk}). Therefore, there are no such curves for this intermediate value of $t$.

The next case is that $t_1$ corresponds to a de-compression of $\mathcal F_t ^ {T_1}$, and $t_2$ corresponds to a maximum of $T_1$. If we choose $t$ just after $t_1$, we see a compressing disk, $D$, for $\mathcal F_t ^ {T_1}$ lying entirely in the interior of some tetrahedron. $\mathcal F_t$ separates $M$ into two components, $A$ and $B$, and suppose $D \subset A$. If the lemma is not immediately true for this value of $t$, then we see either a non-normal arc, or a simple closed curve which is essential in $\mathcal F_t ^ {T_1}$, which cuts off a disk, $D'$ of some face. Theorem \ref{maxwrtk} implies that $D'$ must also be on side $A$, and hence, must be a low disk. As before, we can find a high disk in a face for a surface close to $t_2$, and so we conclude there is some intermediate value where the lemma must be true.

Now we must consider the case when $t_1$ corresponds to a de-compression, and $t_2$ corresponds to a compression. But as above, if the lemma fails to be true for values of $t$ near $t_1$, then we see a low disk in the face of some tetrahedron. Similarly, if the lemma is false for $t$ near $t_2$ we will see a high disk lying in some face. So by the same argument, there is an intermediate value where the lemma is true.

If $\partial M =\emptyset$, then all remaining cases are symmetric to the ones discussed above. If $\partial M \ne \emptyset$, then we need to consider what happens when $t_1$ corresponds to a minimum of $T_1$, and $t_2$ corresponds to a $\partial$-compression of $\mathcal F_t$. This is by far the most difficult case. Choose $t$ just before $t_2$, when we see a $\partial$-compressing disk, $D$, for $\mathcal F_t$, contained entirely in some tetrahedron, $\tau$. Note that $\partial D=\alpha \cup \beta$, where $\mathcal F_t \cap D=\alpha$, and $\partial M \cap D=\beta$. Such a disk is a high disk. Now suppose that $\mathcal F_t \cap T_2$ contains some non-normal arc, or some simple closed curve which is essential on $\mathcal F_t^{T_1}$. As before, this leads us to a high or low disk, $D'$. If $D' \cap D=\emptyset$, then it must be a high disk, also. Furthermore, any low disk in a face of some tetrahedron would be disjoint from $D'$, so there must not be any. We can now repeat the argument given above, to find an intermediate value with no high or low disks in the faces of any tetrahedron.

        \begin{figure}[htbp]
        \psfrag{d}{$\delta$}
        \psfrag{g}{$\gamma$}
        \psfrag{b}{$\beta$}
        \vspace{0 in}
        \begin{center}
        \epsfxsize=5 in
        \epsfbox{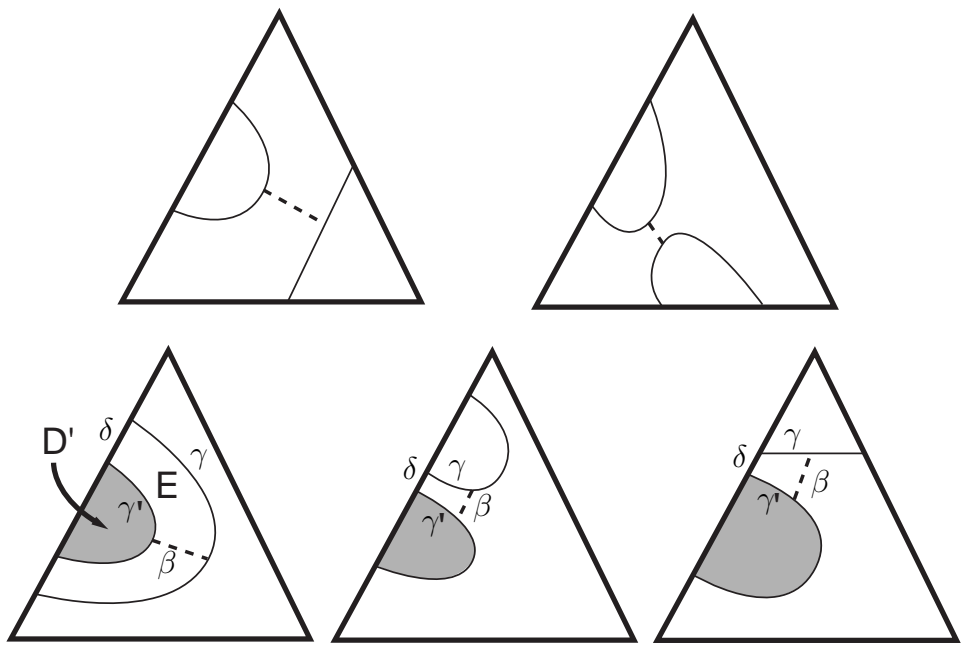}
        \caption{Possibilities for $\mathcal F_t \cap \Delta$, when $t_1$ corresponds to a minimum of $T_1$, and $t_2$ corresponds to a $\partial$-compression.}
        \label{f:boundarysx}
        \end{center}
        \end{figure}

If, however, $D' \cap D \ne \emptyset$, then $D'$ may be a low disk. This leads us to several possibilities. Let $\Delta$ be the face of $\tau$ which contains $\beta$. First of all, if $D'$ is a compressing disk for $\mathcal F_t^{T_1}$, then we see a compressing disk on one side which meets a $\partial$-compressing disk on the other in a single point. These disks can be cancelled, reducing the complexity of $\mathcal F$, and showing that $\mathcal F$ was not locally mini-Lmax. We are left with the possibility that $D'$ is a $\partial$-compressing disk for $\mathcal F_t^{T_1}$. All such configurations are shown in figure \ref{f:boundarysx}. In the bottom three diagrams there is a disk, $E \subset \Delta$, such that $\partial E=\delta \cup \gamma \cup \beta \cup \gamma '$, where $\delta \subset T_1$, $\gamma, \gamma '\subset \mathcal F_t \cap \Delta$, and $D \cap E =\beta$. Note that $E \cup D$ is a $\partial$-compressing disk for $\mathcal F_t^{T_1}$, where $(E \cup D) \cap \mathcal F_t = \gamma \cup \alpha \cup \gamma '$, and $(E \cup D) \cap T_1 = \delta$. If we push $E \cup D$ off of $\Delta$, then we obtain a $\partial$-compressing disk for $\mathcal F_t^{T_1}$ on the opposite side of $\mathcal F_t$ as $D'$, and disjoint from $D'$. This is a contradiction. We conclude that the only possibilities for non-normal arcs are those depicted at the top of figure \ref{f:boundarysx}.

There are still two more cases for $t_1$ and $t_2$, when $t_1$ corresponds to a de-compression or a $\partial$-decompression of $\mathcal F_t$, and $t_2$ corresponds to a $\partial$-compression. These are all similar to those treated above, so they are left as exercises to the reader.
\end{proof}

If we begin with an arbitrary complete collection of maximal leaves for $\mathcal F ^ {T_1}$, then successive applications of Lemma \ref{nlarcs} provides us with a complete collection which intersects every 2-simplex in normal arcs and inessential simple closed curves, with the possible exception of at most 2 non-normal arcs. Suppose $\mathcal F_t$ is a leaf in this collection, and $\gamma$ is an innermost inessential simple closed curve of $\mathcal F_t \cap T_2$. Then $\gamma$ bounds a disk, $D_1$, in the face of some tetrahedron, and a disk, $D_2$, on $\mathcal F_t ^ {T_1}$. Hence, we see a bubble for $\mathcal F_t$. We now invoke Lemma \ref{nobubbles} again to get rid of these bubbles. The result is a foliation in which there is a complete collection of maximal leaves (with respect to $T_1$), where every element of this collection intersects every 2-simplex in normal arcs, and at most 2 non-normal arcs. We shall work with this foliation for the remainder of this section, and we shall assume that $\mathcal F_t$ is an element of our complete collection of maximal leaves.

\begin{lem}
\label{anhalftube}
If $\mathcal F_t \cap T_2$ contains a non-normal arc, then $\mathcal F_t$ is almost normal.
\end{lem}

\begin{figure}[htbp]
        \vspace{0 in}
        \begin{center}
        \epsfxsize=5 in
        \epsfbox{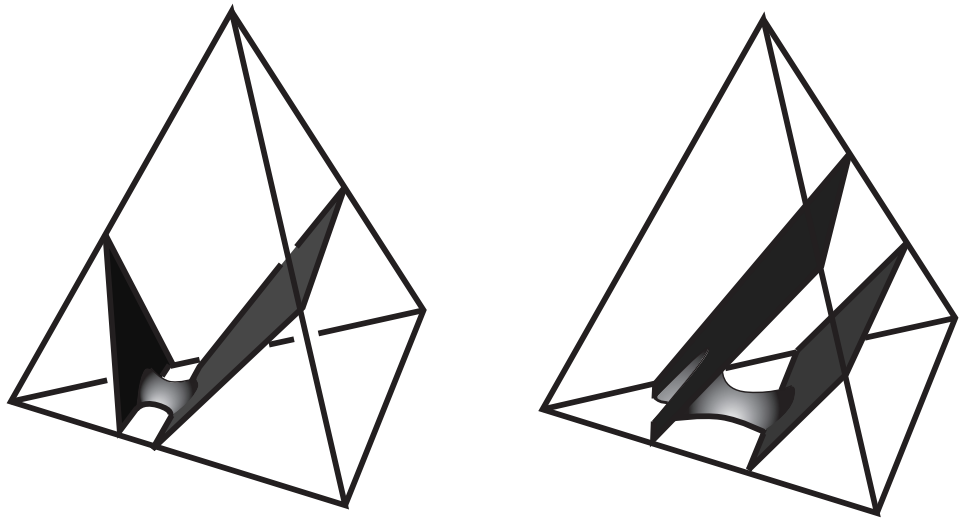}
        \caption{Some possibilities for normal disks, connected by a band which runs along $\partial M$.}
        \label{f:halftube}
        \end{center}
        \end{figure}

\begin{proof}
This situation can only arise when there is a $\partial$-compression of $\mathcal F_t$, as described in the proof of Lemma \ref{nlarcs}. That is, there are two values of $t$, namely $t_1$ and $t_2$, such that $t_1$ somehow corresponds to an increase in $c(\mathcal F^{T_1})$, and $t_2$ corresponds to this $\partial$-compression. (Of course, we may have $t_1$ and $t_2$ switched, but a symmetric argument will hold). Let $t_+$ be some number just after $t_2$. The difference between $\mathcal F_t$ and $\mathcal F_{t_+}$ is that the $\partial$-compression has happened. It is easy to show that if there are any bubbles for $\mathcal F_{t_+}$, then there would be one for $\mathcal F_t$, which there is not. Also, any high or low disk for $\mathcal F_{t_+}$ would be a high or low disk for $\mathcal F_t$, which would be disjoint from the $\partial$-compression. But the boundary compression itself is a high disk, so we cannot see a low disk for $\mathcal F_{t_+}$. Also, as in the proof of Lemma \ref{nlarcs}, there is a low disk for $\mathcal F_t$ which meets the $\partial$-compression in a point. It is also easy to show that any high disk other than the $\partial$-compression would miss this low disk, which is again a contradiction. We conclude that there are no bubbles or high or low disks for $\mathcal F_{t_+}$, and therefore, as in Theorem \ref{normal}, it is normal. Now, $\mathcal F_t$ can be obtained from $\mathcal F_{t_+}$ by un-doing a $\partial$-compression. The picture must be a surface which consists of all normal disks, except for some pair which is connected together by a band that runs along $\partial M$, as in figure \ref{f:halftube}. Such a surface is almost normal.
\end{proof}

\begin{lem}
\label{348}
If $\mathcal F_t \cap T_2$ consists of all normal arcs, then $\mathcal F_t$ meets the boundary of every tetrahedron in normal curves of length 3, 4, and at most one curve on at most one tetrahedron of length 8.
\end{lem}

This lemma is taken straight from \cite{thompson:94}. We refer the reader to this paper for its proof. The necessary assumptions are that $\mathcal F_t$ meets every tetrahedron in normal arcs, and that there are no disjoint high and low disks for $\mathcal F_t$.

\begin{thm}
\label{an}
$\mathcal F_t$ is almost normal.
\end{thm}

\begin{proof}
We now assume that $\mathcal F_t \cap T_2$ is a collection of normal arcs. Let $\tau$ be some tetrahedron in $T$. Let $S$ be a copy of $\partial \tau$, pushed slightly into $\tau$. Now, choose a complete collection of compressing disks for $S \backslash \mathcal F_t$ in $\tau \backslash \mathcal F_t$, and surger $S$ along this collection. We obtain in this way a collection of spheres, $\{ S_1, ... ,S_n \}$. $S_i$ bounds a ball, $B_i$, in $\tau$, and by definition, $\partial B_i \backslash \mathcal F_t$ is incompressible in the complement of $\mathcal F_t ^ {T_1}$ in $M ^ {T_1}$. These are the conditions necessary to apply Theorem 2.1 from \cite{scharlemann:97}. Note that this Theorem is stated only for closed strongly irreducible Heegaard surfaces, but the proof works for Heegaard surfaces with boundary, as in our setting. Hence, there is no problem with the application of this Theorem to $\mathcal F_t ^ {T_1}$. The conclusion is that inside each $B_i$, $\mathcal F_t$ is a connected surface, which looks like the neighborhood of a graph which is the cone on some collection of points in $\partial B_i$. So, in particular, if $\mathcal F_t \cap \partial B_i$ is a single curve, then it bounds a disk in $B_i$, and hence so does the corresponding curve in $\partial \tau$.

        \begin{figure}[htbp]
        \vspace{0 in}
        \begin{center}
        \epsfxsize=5 in
        \epsfbox{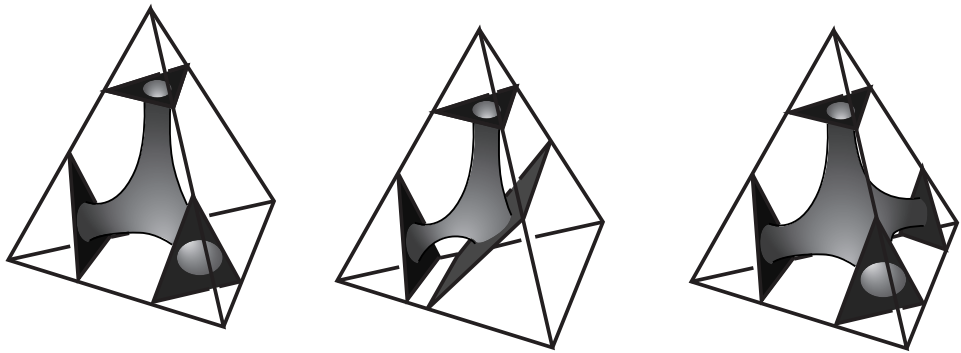}
        \caption{Possibilities when $\mathcal F_t \cap \partial B_i$ consists of 3 or more curves.}
        \label{f:2Tubes}
        \end{center}
        \end{figure}

Suppose there is some $i$ such that $\mathcal F_t \cap \partial B_i$ consists of three or more curves, of length 3 or 4. The only ways this can happen are shown in figure \ref{f:2Tubes}. In all cases we see a compressing disk on one side of $\mathcal F_t$ which is disjoint from a high or low disk on the other side (see figure \ref{f:HLComp}). This cannot happen in a maximal leaf.

        \begin{figure}[htbp]
        \vspace{0 in}
        \begin{center}
        \epsfxsize=3.25 in
        \epsfbox{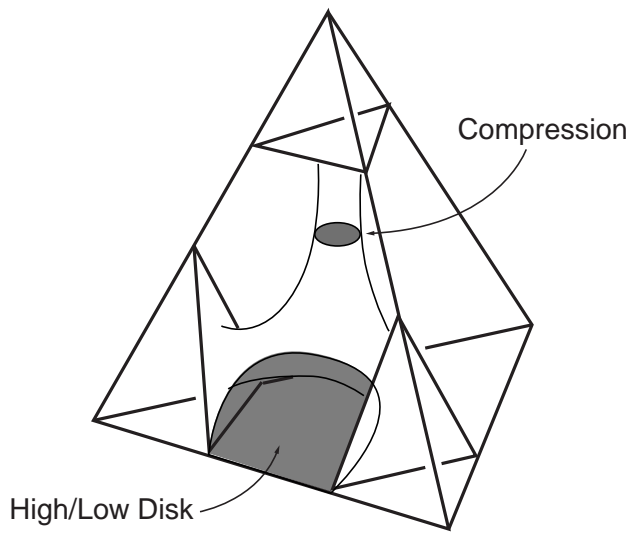}
        \caption{A disjoint compression and High or Low Disk.}
        \label{f:HLComp}
        \end{center}
        \end{figure}

Now suppose that for some $i$, $\mathcal F_t \cap \partial B_i$ consists of two normal curves, of length 3 or 4. \cite{scharlemann:97} tells us that the picture must be two normal disks, tubed together by a single unknotted tube, as in figure \ref{f:1Tube}. Note that in this situation, we see a high or low disk on one side, and a compressing disk on the other. Hence, there cannot be more than one place where we see this picture. Otherwise, we'd see either two disjoint compressing disks on opposite sides, or a compressing disk on one side disjoint from a high or low disk on the other. Neither of these situations can happen for a maximal leaf.

        \begin{figure}[htbp]
        \vspace{0 in}
        \begin{center}
        \epsfxsize=3.75 in
        \epsfbox{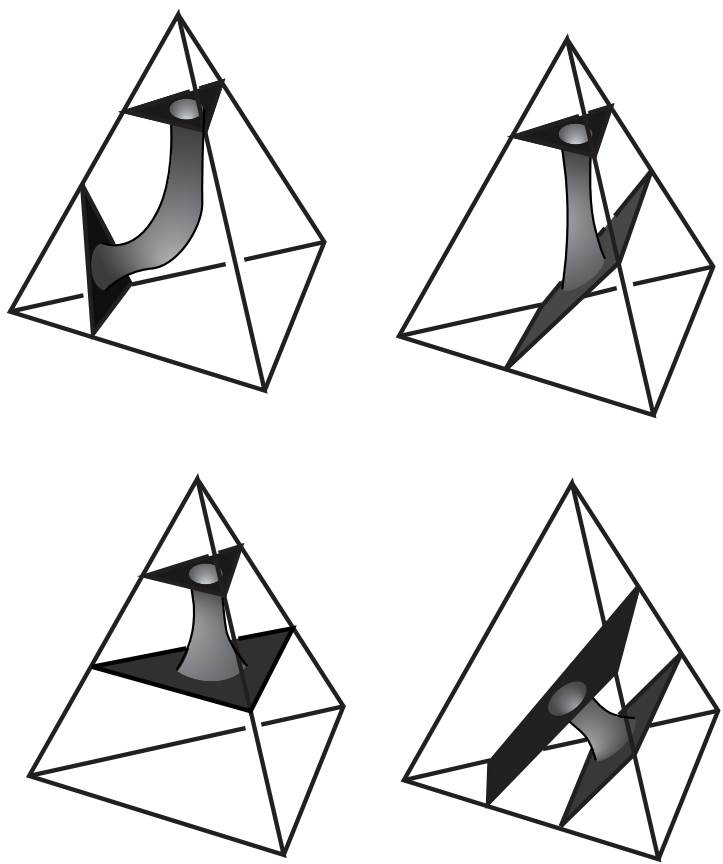}
        \caption{Possibilities when $\mathcal F_t \cap \partial B_i$ consists of 2 curves.}
        \label{f:1Tube}
        \end{center}
        \end{figure}

Furthermore, suppose $\mathcal F_t \cap \partial \tau$ contains a curve of length 8. Then we see a high or low disk on both sides as in figure \ref{f:Octagon}, and hence, there cannot be a tube anywhere else (including attached to this disk!).

        \begin{figure}[htbp]
        \vspace{0 in}
        \begin{center}
        \epsfxsize=3 in
        \epsfbox{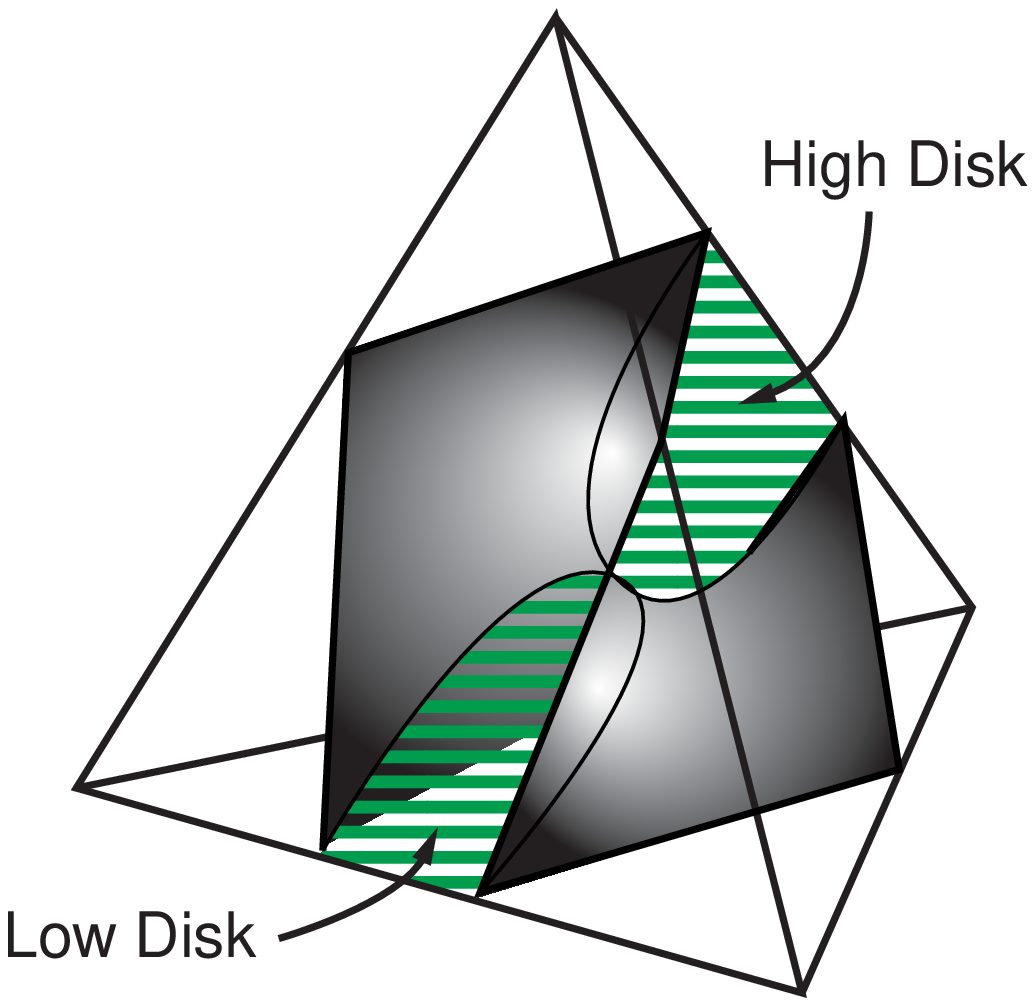}
        \caption{High and Low Disks for an Octagonal piece of $\mathcal F_t$.}
        \label{f:Octagon}
        \end{center}
        \end{figure}

We conclude that $\mathcal F_t$ is made up of all normal disks, with the exception of either a single disk with a boundary curve of length 8, OR a single place where there are two normal disks tubed together by an unknotted tube. This is the precise definition of an almost normal surface.

Our proof is complete by noting that there must be an octagonal disk or a tube {\it somewhere}, because $\mathcal F_t$ is a maximal leaf, and hence there is at least one compression or high or low disk on both sides. If there were no tubes or octagons, then we would not have this.
\end{proof}

As a special case of the Theorem \ref{an}, we obtain a result of Rubinstein \cite{rubinstein:93} and Stocking \cite{stocking:96}, which includes a generalization to $\partial$-Heegaard splittings.

\begin{cor}
\label{hs}
Any (quasi-)strongly irreducible ($\partial$-)Heegaard surface is isotopic to an almost normal surface.
\end{cor}

\begin{proof}
As we have previously seen, any (quasi-)strongly irreducible ($\partial$-)Heegaard surface can be realized as a maximal leaf in a locally mini-Lmax foliation, $\mathcal F$, of $M$. Let $t_1$ and $t_2$ be consecutive critical values, which ``sandwich" the maximal leaf. Hence, at $t_1$ we see $\mathcal F_t$ de-compress, and at $t_2$ we see a compression.

Now, make $T_1$ mini-Lmax with respect to $\mathcal F$. At $t_1$ we still see a de-compression, and so there is still an increase for $c(\mathcal F_t^{T_1})$. Likewise, we still see a compression at $t_2$, and so $c(\mathcal F_t^{T_1})$ still decreases there. Hence, somewhere in between $t_1$ and $t_2$ there is at least one maximal leaf for $\mathcal F ^ {T_1}$. By Theorem \ref{an}, this leaf is an almost normal surface in $M$. But since it is between $t_1$ and $t_2$, it is a maximal leaf for $\mathcal F$, and so it is isotopic to the original (quasi-)strongly irreducible ($\partial$-)Heegaard surface.
\end{proof}

Actually, the full power of Theorems \ref{normal} and \ref{an} lie in the following corollary, which is a strict generalization of the previous result. Recall from \cite{st:94} that a thin decomposition of $M$ is an alternating sequence of incompressible and strongly irreducible surfaces.

\begin{cor}
Any thin decomposition of $M$ can be realized as an alternating sequence of normal and almost normal surfaces.
\end{cor}

\begin{proof}
Any thin decomposition of $M$ is an example of the maximal and minimal leaves of a locally mini-Lmax foliation of $M$. As before, make $T_1$ mini-Lmax with respect to this foliation. By the techniques in the proof of Corollary \ref{hs}, we can easily show that for every minimal (maximal) leaf of $\mathcal F$ there is a minimal (maximal) leaf of $\mathcal F ^ {T_1}$, and hence a normal (almost normal) representative.
\end{proof}

\section{Applications}
\label{s:app}

This section focuses on using the previous results to find normal and almost normal surfaces in knot complements. Our first two Theorems deals with knots which have hyperbolic exteriors. To this end, we will need the following technical Lemma:

\begin{lem}
\label{l:efficient}
If $K$ is a knot in a 3-manifold, $M$, such that $M \backslash K$ admits a complete, hyperbolic structure, then there is a triangulation of $M^K$ in which there is a finite, constructable set of normal and almost normal surfaces of any given Euler characteristic.
\end{lem}

We begin with knots like the ones described in Example \ref{e:thinknot}. That is, suppose $K$ is some knot in $M=S^3$, for which thin position corresponds to bridge position. So, there is some level 2-sphere, which we shall call a {\it bridge sphere}, in $S^3$ which separates all of the maxima from the minima. Define the {\it bridge number of $K$}, $b_K$, to be half of the minimal number of intersections of all possible bridge spheres with $K$. We now apply the results of the previous sections to prove the following theorem:

\begin{thm}
\label{thinalg}
If $K$ is a hyperbolic knot, then there is an algorithm which will either determine the bridge number of $K$, or determine that there is a closed incompressible surface in the complement of $K$.
\end{thm}

\begin{proof}
By Thompson's theorem \cite{thompson:97} we know that if $M^K$ does not contain a meridional, planar, incompressible surface, with fewer boundary components than the width of {\it any} presentation of $K$, then $K$ has a thin presentation which is also bridge. This is a condition we can algorithmically check by \cite{jr:98}, since such a planar incompressible surface has bounded Euler characteristic. If $M$ contains such a surface, then $M$ also contains a closed incompressible surface (see \cite{thompson:97}). So we may now proceed assuming that $K$ has a thin presentation which is also bridge, and show that if this is the case, one can always determine the bridge number of $K$.

Triangulate the complement of $K$ in $S^3$, so that $T_0 \subset \partial M^K$, and so that there are no normal 2-spheres in $M^K$ or non-boundary parallel normal tori. Since $M^K$ is hyperbolic, such a triangulation exists by a result of Casson (see \cite{lackenby:99} for a proof). By the remarks in Example \ref{e:thinknot}, we know that there is a bridge sphere, $S \subset M$, which realizes the minimal number of intersections with $K$, such that $S^K$ is a strongly irreducible $\partial$-Heegaard splitting for $M^K$. We now apply Corollary \ref{hs} to make $S^K$ almost normal.

The algorithm proceeds as follows: First, given any picture of $K$, we can compute $b$, an upper bound for $b_K$, by counting the number of maxima in the picture. Since the set of normal surfaces are finitely generated, Euler characteristic is additive, and there are no normal 2-spheres or non-boundary parallel tori, it follows that there is a finite, constructable set of almost normal, meridional, planar surfaces in $M^K$ with at most $b$ boundary components. We can now look at each, and decide whether or not it is a punctured bridge sphere, by checking to see of it compresses completely to both sides. Among all planar surfaces that do, choose one, $S$, with fewest number of boundary components. This will be a punctured bridge sphere for $K$, which realizes the minimal number of intersections with $K$. $b_K$ then equals half the number of boundary components of $S$.
\end{proof}

{\bf Technical Note.} The result from \cite{jr:98} which we use here says that given a manifold with one boundary component, with no essential 2-spheres, disks, tori or annuli, then there is a triangulation in which all summands with non-negative Euler characteristic of arbitrary normal surfaces can be ignored. Since Euler characteristic is additive when adding normal surfaces, we see that there are a finite number of normal and almost normal surfaces of bounded Euler characteristic. It is likely that similar results hold for manifolds with essential tori and annuli. In this case, we would be able to remove the assumption of hyperbolicity from Theorem \ref{thinalg}.

\bigskip

We can generalize Theorem \ref{thinalg} to knots in manifolds other than $S^3$. Suppose $K$ is some knot embedded in an orientable, irreducible, non-Haken 3-manifold, $M$. The {\it bridge number of $K$}, $b_K$, is the minimal number of maxima of $K$, among all embeddings of $K$ which are in bridge position with respect to any minimal genus strongly irreducible Heegaard splitting of $M$ (see Theorem \ref{t:generalthompson}).

\begin{thm}
Let $M$ be a closed, orientable, irreducible, non-Haken 3-manifold, and let $K$ be a knot in $M$ with hyperbolic exterior, which is not isotopic onto any minimal genus strongly irreducible Heegaard splitting. Then there is an algorithm which will either determine the bridge number of $K$, or find a meridional, incompressible, $\partial$-incompressible surface in the complement of $K$, which has genus less than or equal to that of $M$.
\end{thm}

\begin{proof}
Let $H$ be some minimal genus strongly irreducible Heegaard splitting of $M$. The proof follows exactly that of the previous Theorem, where we substitute Theorem \ref{t:generalthompson} for \cite{thompson:97}. Theorem \ref{t:generalthompson} tells us that if $M^K$ does not contain a meridional incompressible surface, with fewer boundary components than the width of {\it any} embedding of $K$, and genus smaller than or equal to $H$, then there is an embedding of $K$ which is both locally mini-Lmax, and bridge with respect to $H$. This is a condition we can algorithmically check by \cite{jr:98}, since such an incompressible surface has bounded Euler characteristic. We now proceed assuming that $K$ has an embedding which is both locally mini-Lmax and bridge with respect to $H$, and show that if this is the case, one can always determine $b_K$.

Triangulate the complement of $K$ in $M$, so that $T_0 \subset \partial M^K$, and so that there are no normal 2-spheres in $M^K$. Theorem \ref{maxwrtk} implies that if $K$ is embedded so that it is both locally mini-Lmax and bridge with respect to $H$, then $H^K$ is a quasi-strongly irreducible $\partial$-Heegaard splitting for $M^K$. We now apply Corollary \ref{hs} to make $H^K$ almost normal.

The algorithm proceeds as follows: First, use \cite{rubinstein:95} to determine the Heegaard genus of $M$. Then, given any embedding of $K$, we can compute $b$, an upper bound for $b_K$, by counting the number of maxima in the picture. By \cite{jr:98}, there is a finite, constructable set of almost normal, meridional, surfaces in $M^K$ with at most $b$ boundary components, and genus equal to that of $M$. We can now look at each, and decide two things: first, whether or not it compresses completely to both sides in $M^K$, and second, whether or not the corresponding surface divides $M$ into two handle-bodies. Among all surfaces that satisfy both, choose one, $H'$, with fewest number of boundary components. This will be a punctured minimal genus Heegaard splitting of $M$, which realizes the minimal number of intersections with $K$. $b_K$ then equals half the number of boundary components of $H'$.
\end{proof}

Before proceeding to the next theorem, we need a new definition.

\begin{dfn}
$\Sigma$ is an {\it untelescoped Heegaard decomposition} of $M$ if $\Sigma$ is the disjoint union of maximal leaves in a locally mini-Lmax foliation of $M$.
\end{dfn}

\begin{thm}
\label{exteriors}
Let $M$ be an irreducible 3 manifold, and $K$ a knot in $M$. Let $M^K$ denote $M$ with a regular neighborhood of $K$ removed. Then one of the following is true:
    \begin{itemize}
    \item $M^K$ contains a meridional almost normal surface.
    \item $M^K$ contains a meridional normal surface, which is planar, incompressible, and $\partial$-incompressible.
    \item $M^K$ contains an essential normal 2-sphere.
    \item $K$ is isotopic onto every untelescoped Heegaard decomposition of $M$.
    \end{itemize}
\end{thm}

\begin{proof}
First, if $K$ is locally untangled and cannot be isotoped onto a leaf of some mini-Lmax foliation, $\mathcal F$, then Theorems \ref{maxwrtk} and \ref{an} say that the maximal leaves of $\mathcal F^K$ can be realized as a union of almost normal surfaces. Since $K$ must have a minimum and a maximum with respect to the height function which induces $\mathcal F$, there must be a maximal leaf of $\mathcal F^K$ which hits $K$. Hence, if $K$ is locally untangled and cannot be isotoped onto a leaf of $\mathcal F$, we have an almost normal meridional surface.

If $K$ is locally tangled, then by definition there is a ball, $B \subset M$, such that $\partial B^K$ is incompressible and $\partial$-incompressible in $M^K$, or such that $K \subset B$. In either case, $\partial B$ can be made normal.

Lastly, we have the possibility that $K$ can be isotoped onto a leaf of every foliation, $\mathcal F$. This is equivalent to saying that $K$ can be isotoped onto any untelescoped Heegaard decomposition of $M$.
\end{proof}

For the remainder of the paper, let $X$ denote an irreducible, orientable 3-manifold such that $\partial X$ consists of a single torus. A {\it slope}, $\alpha$, is an isotopy class of essential simple closed curves on $\partial X$. By a {\it Dehn filling along $\alpha$}, we mean the manifold, $X(\alpha )$, obtained from $X$ by gluing a solid torus, $T$, to $\partial X$, in such a way that $\alpha$ bounds a disk in $T$. Finally, let $K$ denote the core of $T$ in $X(\alpha )$.

We would now like to cite the recent work of Jaco and Sedgwick \cite{js:98}:

\begin{thm}
\label{fin}
(Jaco-Sedgwick) If $X$ is triangulated in such a way so as to induce a triangulation of $\partial X$ with one vertex, then there is a finite, constructable set of normal curves on $\partial X$ that can be the boundary of a normal or almost normal surface.
\end{thm}

Unfortunately, some of the slopes that can be the boundaries of almost normal surfaces in this paper are not normal curves. This happens precisely when the exceptional disk is of the type depicted in Figure \ref{f:halftube}. It is therefore necessary to prove the following:

\begin{thm}
\label{finite}
If $X$ is triangulated in such a way so as to induce a triangulation of $\partial X$ with one vertex, then there are a finite number of slopes on $\partial X$ that can be the boundary of a normal or almost normal surface.
\end{thm}

\begin{proof}
The proof follows directly from Theorem \ref{fin}. As noted above, the only time where Theorem \ref{fin} is not sufficient is if we have an almost normal surface with an exceptional piece of the type depicted in Figure \ref{f:halftube}, so that its boundary is not a normal curve. Such a curve is made of all normal arcs, except for in one place, where we see either of the two pictures at the top of Figure \ref{f:boundarysx}. We will call such a curve {\it almost normal}. To prove the Theorem, it suffices to show that every almost normal curve that can be the boundary of an almost normal surface can be obtained from some normal boundary slope in one of a finite number of possible ways.

Note that the exceptional pieces depicted in Figure \ref{f:halftube} each $\partial$-compress to a pair of normal disks. Hence, an almost normal surface which contains one of these exceptional pieces will $\partial$-compress to a normal surface. On $\partial X$, the $\partial$-compression looks like a band sum along the dashed line in Figure \ref{f:boundarysx}. The dual process is to take a normal curve, and band sum along an arc which connects two different normal arcs. But for each normal curve, there are a finite number of pairs of normal arcs. The result now follows from Theorem \ref{fin}.
\end{proof}

This result, together with Theorem \ref{exteriors}, immediately gives us:

\begin{cor}
\label{Dehn}
If $X$ is irreducible, then for all but finitely many slopes, $\alpha$, on $\partial X$, $K$ can be isotoped onto every untelescoped Heegaard decomposition of $X(\alpha )$.
\end{cor}

This corollary is closely related to a Theorem of Rieck \cite{rieck:97}, who proves that for all but finitely many fillings, $K$ can be isotoped onto strongly irreducible Heegaard splittings of bounded genus. Our result is stronger in the sense that we have removed the assumption of bounded genus, but weaker, because we do not find an explicit bound on the number of fillings where $K$ cannot be isotoped onto a leaf of some foliation. This result is also closely related to Theorem 0.1 of \cite{mr}.

Our work also gives short proofs of the following two Corollaries, which have also been proved using different techniques by Jaco and Sedgwick \cite{js:98}:

\begin{cor}
\label{c:knotrecognition}
There is an algorithm to find any slope, $\alpha$, such that $X(\alpha )$ is homeomorphic to $S^3$.
\end{cor}

\begin{proof}
Suppose $\alpha$ is a slope on $\partial X$ such that $X(\alpha )$ is homeomorphic to $S^3$. By \cite{thompson:97}, either $X$ contains a planar, incompressible, $\partial$-incompressible surface with boundary slope $\alpha$, or thin position for $K$ is the same as bridge position. In the former case, there will be a normal surface in $X$ with boundary slope $\alpha$. In the latter case, a bridge 2-sphere for $K$ will be a strongly irreducible $\partial$-Heegaard splitting of $X$ (see Example \ref{e:thinknot}), which can be made almost normal by Corollary \ref{hs}. In either case, we get a normal or almost normal surface in $X$, with boundary slope $\alpha$. We now apply Theorem \ref{finite}, which says there is a finite, constructable set of such slopes. For each such slope, we can form $X(\alpha )$, and decide whether or not the manifold is $S^3$ by \cite{rubinstein:95}.
\end{proof}

\begin{cor}
There is an algorithm to find all lens space fillings of $X$, or any filling homeomorphic to $S^2 \times S^1$.
\end{cor}

\begin{proof}
Suppose $\alpha$ is a slope on $\partial X$ such that the Heegaard genus of $X(\alpha )=1$. Our first case to consider is when $K \subset B$, for some ball $B \subset X(\alpha )$. Then a prime decomposition for $X$ will be the connect sum of a manifold with Heegaard genus 1, and the complement of a knot in $S^3$. We can recognize the former by \cite{rubinstein:93}, and the latter by Corollary \ref{c:knotrecognition}.

If $X$ is irreducible, then Theorem \ref{exteriors} says that either $K$ can be isotoped to lie on the Heegaard torus in $X(\alpha )$, or $X$ contains a normal or almost normal surface with boundary slope $\alpha$. In the latter case, Theorem \ref{finite} says that there are a finite number of possibilities for $\alpha$, and by \cite{rubinstein:93}, we can recognize which ones of these correspond to Dehn fillings which have Heegaard genus 1.

If, on the other hand, $K$ can be isotoped to lie on a Heegaard
torus for $X(\alpha )$, then $X$ contains an essential annulus,
which can be normalised. Furthermore, cutting along this annulus
yields two solid tori. By \cite{jt:95}, we can decide whether or
not $X$ contains an essential annulus, and by \cite{haken:61}, we
can tell if an irreducible manifold is a solid torus. Hence, we can decide if
$X$ admits such a decomposition a priori.
\end{proof}

\appendix

\section{Proof of Theorem \ref{t:bcg}}
\label{a:bcg}

First, let $M$ be any 3-manifold with a $\partial$-Heegaard splitting, $W \cup _F W'$. Note that $M$ can be described as follows: Begin with $F \times I$, and attach 2-handles and half 2-handles to $F \times \partial I$. Finally, cap off 2-sphere boundary components on each side with 3-balls, and add any disk components of the boundary to $\partial _0 M$. Then $W$ and $W'$ are the submanifolds obtained by cutting $M$ along $F \times \{1/2\}$. Also note that if $F \times \{1/2\}$ is a quasi-strongly irreducible $\partial$-Heegaard splitting surface, then so is $F \times \{t\}=F_t$, for any $t \in (0,1)$. If we cut $M$ along $F_t$, we obtain two $\partial$-compression bodies, which we shall denote $W_t$ and $W'_t$ (where $W_t$ is the one which contains $F_0$).

Now, suppose $D$ is a compressing or $\partial _0$-compressing disk for $\partial _-M$. Let {\bf D, D'} be complete collections of disjoint 2-handles and half 2-handles attached to $F_0$ and $F_1$, respectively (in the sense that $F_0$ compressed and $\partial$-compressed along the cores of all the handles in {\bf D} is a 2-sphere, a disk parallel to $\partial _0 M$, or a surface parallel to $\partial _- M$).

Let $D \cap \partial _0 M=\partial _0 D$. Let $\pi _I:F \times I \rightarrow I$ denote the projection map. Let $m=|int(\partial _0 D) \cap (F \times \partial I)|$, and let $n$ equal the number of critical point of $\partial _0 D$ with respect to $\pi _I$. We now assume that $D$ was chosen so that $(m,n)$ is minimal. Note that if $D$ is a compressing disk (as opposed to a $\partial$-compressing disk) for $\partial _- M$, then $(m,n)=(0,0)$. Now, after isotopies, compressions, and $\partial$-compressions of $D$, we may also assume that each component of $D \cap (\bf D \cup {\bf D'})$ is a disk which lies in some element of {\bf D} or {\bf D'}, and is parallel to its core. Such a move can only lower $m$, so there is no problem in continuing in our assumption that $(m,n)$ is minimal. Note that this puts $D$ into a position where $D \cap \partial _-M=\partial _-D$ lies on $F_0$ (say), and misses all of the regions where the handles of {\bf D} are attached.

\begin{clm}
\label{c:n=0}
n=0.
\end{clm}

\begin{proof}
If $D$ is a compressing disk for $\partial _- M$, there is nothing to prove. So we begin by assuming that $D$ is a $\partial _0$-compressing disk, and $m=0$. That is, $\partial _0 D$ lies entirely in $\partial F \times I$. Since both endpoints of $\partial _0 D$ must lie on the same component of $F \times \partial I$, it must be that $\partial _0 D$ co-bounds a subdisk, $E$, of $\partial F \times I$. Now, we can use $E$ to isotope $D$ so that it becomes a compressing disk for $\partial _- M$. This shows that $(m,n)$ was not minimal for our original choice of $D$.

If $m>0$, then let $\alpha$ be a component of $\partial _0 D \cap (\partial F \times I)$ which contains a critical point of $\partial _0 D$ with respect to $\pi _I$. If both endpoints of $\alpha$ lie on the same component of $F \times \partial I$, then again $\alpha$ co-bounds a subdisk, $E$, of $\partial F \times I$. Let $\alpha '$ be an outermost arc of $\partial _0 D \cap E$, and let $E'$ be the subdisk of $E$ which it bounds. Then we can use $E'$ to isotope $D$, lowering $m$ by two.

If, on the other hand, the endpoints of $\alpha$ lie on different components of $F \times \partial I$, then $\alpha$ can be straightened to an arc of the form $p \times I$, where $p$ is some point of $\partial F$. This lowers $n$, contradicting our original assumption of minimality.
\end{proof}

Let $t \in (0,1)$, and suppose $\gamma$ is an arc of $F_t \cap D$ which is outermost on $D$. $\gamma$ cuts off a subdisk, $D'$, of $D$.

\begin{clm}
\label{c:honest}
If $D'$ does not contain any simple closed curves of $F_t \cap D$ which are essential on $F_t$, then $D'$ is isotopic to an honest $\partial _0$-compressing disk for $F_t$.
\end{clm}

\begin{proof}
By an innermost disk argument, we can isotope $D'$ to remove simple closed curve components of $F_t \cap D'$, which are inessential on $F_t$. After doing this, $D'$ is entirely contained in $W_t$ (say). Let $\partial _0 D'=D' \cap \partial _0 W_t$. We now claim that $\partial _0 D'$ is essential on $\partial _0 W_t$. Suppose not. Then there exists a disk, $E \subset \partial _0 W_t$, such that $\partial E =\partial _0 D' \cup \alpha$, where $\alpha=E \cap F_t$. If $E$ is entirely contained in $F \times I$, then $\partial _0 E =\partial _0 D'$ must contain a critical point with respect to $\pi _I$. This implies that $n>0$, contradicting Claim \ref{c:n=0}. Otherwise, let $\beta$ be an arc of $E \cap (F \times \partial I)$, which is outermost on $E$. Let $E'$ be the subdisk of $E$ cut off by $\beta$. Now, we can use $E'$ to guide an isotopy of $D$ which lowers $m$.

Now, suppose $D' \cap F_t$ is an inessential arc on $F_t$.  Then $D'$ can be isotoped off of $F_t$, to become a compressing disk for $\partial _0 W_t$. Since $\partial _0 W$ is incompressible in any $\partial$-compression body, $W$, this is a contradiction. We conclude that $D'$ must be an honest $\partial _0$-compressing disk for $F_t$ in $W_t$
\end{proof}

Let $\Gamma _t$ be the subcollection of 1-manifolds of $D \cap F_t$ which are essential on $F_t$. An element, $\gamma$, of $\Gamma _t$ is an {\it $H$-curve} if it cuts off a subdisk, $D'$, of $D$, such that $D'$ contains no other element of $\Gamma _t$, and such that a collar of $\gamma$ in $D'$ lies in $W_t$. We define an {\it $L$-curve} similarly, the only difference being that for an $L$-curve, a collar of $\gamma$ in $D'$ must lie in $W'_t$. Note that if an $H$- or an $L$-curve is closed, then it is an innermost loop of $\Gamma _t$ on $D$. If it is an arc, then it is an outermost arc.

By a standard innermost disk/outermost arc argument, and Claim \ref{c:honest}, we can show that any $H$-curve bounds a compressing/honest $\partial _0$-compressing disk for $F_t$ in $W_t$, and any $L$-curve bounds a compressing/honest $\partial _0$-compressing disk for $F_t$ in $W'_t$. Hence, it follows from the quasi-strong irreducibility of $W_t \cup _{F_t} W'_t$ that $F_t$ cannot contain both an $H$-curve and an $L$-curve.

Now, for small $\epsilon$, we know that $F_\epsilon$ contains an $H$-curve. This follows from the fact that $D \cap${\bf D} must be non-empty. Otherwise, $D$ would be disjoint from the core of every handle of {\bf D}. Since $D$ and {\bf D} are on opposite sides of $F_0$, this leads to disjoint compressing/honest $\partial _0$-compressing disks for $F_\epsilon$, contradicting quasi-strong irreducibility. In addition, we know that $F_{1-\epsilon}$ must contain an $L$-curve. Otherwise, $D \cap${\bf D'} would be empty, and $D$ would be a compressing/honest $\partial _0$-compressing disk for $\partial _-W_{1-\epsilon}$ in $W_{1-\epsilon}$ (a contradiction).

We now claim that there exists an interval, $(t_0,t_1)\subset I$, such that for every $t \in (t_0,t_1)$, $F_t$ contains no $H$- or $L$-curves. Note that as $t$ varies, the collection $\Gamma _t$ can only change at saddle tangencies of $D\cap F_t$ (center tangencies only create/destroy curves which are inessential on $F_t$). However, the curves of $D \cap F_t$ just before a saddle tangency can be made disjoint from the curves afterwards. Hence, if there is an $H$-curve before a saddle tangency, there cannot be an $L$-curve afterwards. We conclude that as $t$ varies from $\epsilon$ to $1-\epsilon$, there cannot be an instantaneous transition from $H$-curves to $L$-curves. So there must be an open interval where there are neither.

Now let $t\in (t_0,t_1)$. The fact that there are no $H$- or $L$-curves for $F_t$ immediately implies that $\Gamma _t$ must be empty. Hence, every curve of $D \cap F_t$ must be inessential on $F_t$. We can now apply an innermost disk/outermost arc argument to isotope $D$ so that $D \cap F_t=\emptyset$. This makes $D$ a compressing/honest $\partial _0$-compressing disk for $\partial _-W_t$ in $W_t$, a contradiction.

\bibliographystyle{plain}

\end{document}


\begin{abstract}
We present a complete proof of Theorem 6.3 from ``Heegaard splittings with boundary and almost normal surfaces [Topology and Appl. 116(2001) 153-184]." We also take this opportunity to correct a few minor errors that appeared in that paper. 
\end{abstract}
\maketitle

\pagestyle{myheadings}
\markboth{DAVID BACHMAN}{HEEGAARD SPLITTINGS WITH BOUNDARY (ERRATUM)}

We assume the reader is familiar with \cite{bachman:98}. The proof of Theorem 6.3 of that paper contains the phrase ``Note that $\partial B$ cannot compress away to nothing outside $B$, since $M$ is not homeomorphic to $S^3$" (page 167, lines 4 and 5). While, strictly speaking, this statement is true it certainly requires further proof. This is the main content of the present paper. 

One additional point which we address here is that Theorem 6.3 is stated only for 1-manifolds, but is used later for 1-vertex graphs. Furthermore, some of this graph may lie on $\partial M$, a point which was not emphasized in the original. The introductory paragraphs to Section 8 of \cite{bachman:98} were meant to deal with these technicalities, but some readers may find them somewhat unsatisfying. We will be much more precise with these issues here.

Before getting started the author would also like to apologize for not referencing the work of Yoav Reick and Eric Sedgwick \cite{sr:01} in \cite{bachman:98}, since they had concurrently obtained many of the same results. The author is also grateful to Eric Sedgwick and Saul Schleimer for helpful conversations during the preparation of this note.  

We begin by defining exactly what sort of graphs we will be dealing with.

\begin{dfn}
\rm A graph $K$ embedded in a 3-manifold $M$ is {\it locally untangled} if 
\begin{enumerate}
    \item $K$ has no vertices of valence zero, one, or two, 
    \item the interior of each edge is either contained in $\partial M$ or is disjoint from it, and
    \item for each ball $B \subset M$ with incompressible frontier in $M^K$ there is a point $p \in K$ such that for all sufficiently small $\epsilon$ the frontier of $B$ in $M^K$ is parallel to the frontier of $N_{\epsilon}(p)$ in $M^K$.
\end{enumerate}
\end{dfn}

It follows that a connected locally untangled graph has a single vertex. If not, then let $e$ be an edge connecting distinct vertices. A neighborhood of $e$ in $M$ will be a ball violating the third condition of the definition.

Examples of locally untangled graphs include 1-vertex spines of handlebodies bounded by strongly irreducible Heegaard splittings \cite{scharlemann:97} and 1-skeleta of {\it 0-efficient} triangulations \cite{jr:98}. It is the latter application that we make use of in \cite{bachman:98}.

We now define the complexity of a surface with respect to a locally untangled graph. If $F$ is a connected, properly embedded surface in $M$ and $K$ is a locally untangled graph then we define $c_K(F)$ to be $|F \cap K|$ plus
\begin{enumerate}
    \item $0$ if $F$ is a disk or sphere
    \item $1-\chi(F)$ if $F$ is closed
    \item $\frac{1}{2}-\chi(F)$ otherwise
\end{enumerate}
If $F$ is not connected then we define $c_K(F)$ to be $\sum _i c_K(F_i)$, where the sum ranges over all components $F_i$ of $F$. {\it Maximal} and {\it minimal} leaves (with respect to $K$) of a singular foliation $\mathcal F$ are now defined in the obvious way, as well as the complexity $Lmax_K(\mathcal F)$. 

\begin{dfn} 
\rm Let $\mathcal F$ be a singular foliation arising from a height function $h:M \rightarrow I$. A locally untangled graph $K$ with vertex $v$ is in {\it good position} if $v$ is a local extremum of $h|K$. The graph $K$ is {\it mini-Lmax with respect to $\mathcal F$} if it is in good position, and $K$ cannot be isotoped to reduce $Lmax_K(\mathcal F)$.
\end{dfn}

\begin{figure}[htbp]
        \psfrag{P}{$P$}
        \psfrag{Q}{$Q$}
        \psfrag{M}{$\partial_0 M$}
        \psfrag{+}{$\partial _+ W^K$}
        \psfrag{-}{$\partial _- W^K$}
        \psfrag{0}{$\partial _0 W^K$}
        \vspace{0 in}
        \begin{center}
        \includegraphics[width=4 in
        ]{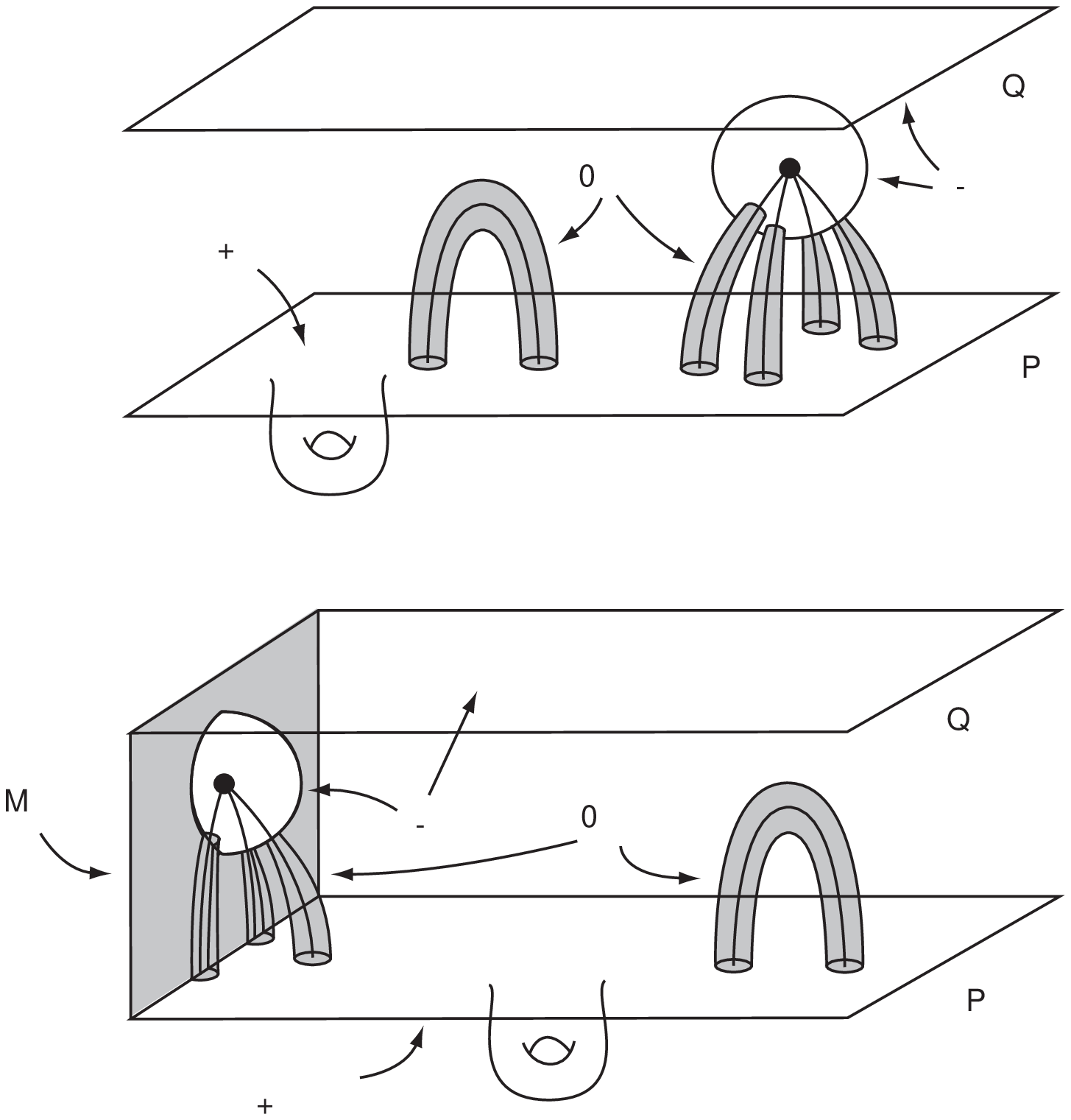}
        \caption{}
        \label{f:bdycompbod}
        \end{center}
        \end{figure}

If $K$ is in good position, $P^{K}$ is a maximal leaf of $\mathcal F^{K}$, and $Q^{K}$ is the next minimal leaf then the region of $M^K$ between $P^K$ and $Q^K$ is a $\partial$-compression body $W^K$. One ambiguity here is that for each $\partial$-compression body there are always many ways to assign $\partial _+$, $\partial _-$, and $\partial _0$. For the region between $P^K$ and $Q^K$ we will always make these assignments in the following way. If the vertex $v$ of $K$ is in $W$ then let $B$ be a neighborhood of $v$ in $M$, and let $S$ denote the frontier of $B^K$ in $M^K$. Now, let 
\begin{enumerate}
    \item $\partial _+ W^K=P^K$, 
    \item $\partial _-W^K=Q^K \cup S$, 
    \item and $\partial _0 W^K$ denote the remaining boundary of $W$.
\end{enumerate}
It follows that $\partial _0 W^K  \subset \partial N(K) \cup \partial_0 M$. If $M$ is closed then these assignments are pictured schematically at the top of Figure \ref{f:bdycompbod}. When $M$ has boundary they are pictured at the bottom. Note that with these assignments a $\partial _0$-compression can run over maxima and minima of $K$, but not over the vertex.

One final correction is to the definition of a locally mini-Lmax foliation. On page 164, line 16 of \cite{bachman:98} we say, ``we shall even refer to $\mathcal F$ as locally mini-Lmax if the maximal leaves are only quasi-strongly irreducible Heegaard surfaces." This sentence should be omitted. That is, a foliation is locally mini-Lmax only if it satisfies the conclusion of Theorem 5.2 of \cite{bachman:98}. 

The correct statement of Theorem 6.3 is as follows. 

\medskip
\noindent
{\bf Theorem 6.3.} {\it Let $M$ be an irreducible 3-manifold other than $S^3$. Let $K$ be a locally untangled graph in $M$ with at least one edge whose interior is not on $\partial M$. Let $\mathcal F_*$ be a locally mini-Lmax foliation of $M$ such that no loop component of $K$ can be isotoped onto a leaf of $\mathcal F_*$. Then there is a locally mini-Lmax foliation $\mathcal F$ with the same maximal and minimal leaves as $\mathcal F_*$ (and hence the same $Lmax$ complexity) such that when $K$ is isotoped to be mini-Lmax with respect to $\mathcal F$ the foliation $\mathcal F^K$ has the following property. If $K$ has a vertex then the maximal leaves of $\mathcal F^K$ are quasi-strongly irreducible $\partial$-Heegaard surfaces for the submanifolds of $M^K$ that arise when we cut along the minimal leaves. If $K$ has no vertex ({\it i.e.} it is a 1-manifold) then the maximal leaves are strongly irreducible, and hence $\mathcal F^K$ is a locally mini-Lmax foliation of $M^K$.}
\medskip

Before proceeding with the proof we will need the following lemmas.

\begin{lem}
\label{l:thinball}
Let $L$ be a leaf of $\mathcal F$ and $S$ be a sphere in $M$ such that $S \cap K \subset L$ and $S \cap L$ is a connected surface. Then there exists a $K' \subset K$ such that $\partial K' \subset S \cap L$ and such that there is an isotopy of $K$ which fixes $K \backslash K'$, and sends $K'$ into $S \cap L$. 
\end{lem}

\begin{proof}
Compress $S$ as much as possible in $M^K$, to obtain a collection of spheres, $R$, in $M$. If any component of $R$ bounds a ball which contains a single unknotted arc of $K$, then there is an innermost such, $R'$. The situation is now similar to the proof of the main result of \cite{thompson:97}. Inside $R'$ there is a disk, $T$, such that $\partial T=\delta \cup \gamma$, where $\delta \subset \partial N(K)$ and $\gamma \subset R'$. We now reverse the compressions used to obtain $R$ from $S$. Each time a compression is reversed, we attach a tube to some components of $R$. These tubes may intersect $T$, but only in its interior. The arc $\delta$ is parallel in $N(K)$ to an arc $K' \subset K$. We can now use $T$ to guide an isotopy of $K$ which sends $K'$ into $S$. As $S$ is a sphere and $S \cap L$ is connected we may do a further isotopy of $K'$ in $S$ to bring $K'$ into $S \cap L$.

Since $K$ is locally untangled, the other possibility is that every component of $R$ is parallel in $M^K$ to the boundary of a neighborhood of the vertex of $K$. This case is similar to the main argument of \cite{thinball}, and the author is greatly indebted to Saul Schleimer for his collaboration on that work.

To fix notation, let $\{D_i\}_{i = 1}^n$ be a sequence of disks, and $\{S_i\}_{i = 0}^n$ the sequence of surfaces, so that $S_0 = S$, $S_n = R$, and $S_i$ is obtained from $S_{i-1}$ by compressing along $D_i$ in the complement of $K$.  That is, remove a small neighborhood, $A_i$, of $\partial D_i$ from $S_{i-1}$. Construct $S_i$ by gluing two parallel copies of $D_i$ onto $\partial A_i$. Denote these by $B_i$ and $C_i$. So $A_n \cup B_n \cup C_n$ bounds a ball homeomorphic to $D_n \times I$. Finally, let $\alpha$ be the image of $\{pt\} \times I \subset D_n  \times I$ by such a homeomorphism. The arc $\alpha$ is a {\em co-core} for the compression $D_n$.

Let $V$ be the component of $S_{n-1}$ which meets $D_n$.  Let $\{R_j\}_{j = 0}^{n}$ denote the components of $R = S_n$, numbered consecutively, so that $R_0$ is farthest from the vertex of $K$. Then there is a $j$ such that compressing $V$ along $D_n$ yields $R_j$ and $R_{j+1}$.  Let $N$ denote the submanifold of $M$ bounded by $R_j$ and $R_{j+1}$. Note $\alpha \subset N$.

Choose a homeomorphism $h:S^2 \times I \rightarrow N$ such that $K \cap N$ is a collection of {\em straight arcs} and $\pi(h^{-1}(B_n)) = \pi(h^{-1}(C_n))$. Here a straight arc is one of the form $h(\{pt\} \times I)$ and $\pi:S^2 \times I \rightarrow S^2$ is projection onto the first factor.

        \begin{figure}[htbp]
        \psfrag{R}{$f_1(U)$}
        \psfrag{K}{$K'$}
        \psfrag{J}{$g$}
        \vspace{0 in}
        \begin{center}
        \includegraphics[width=3 in]{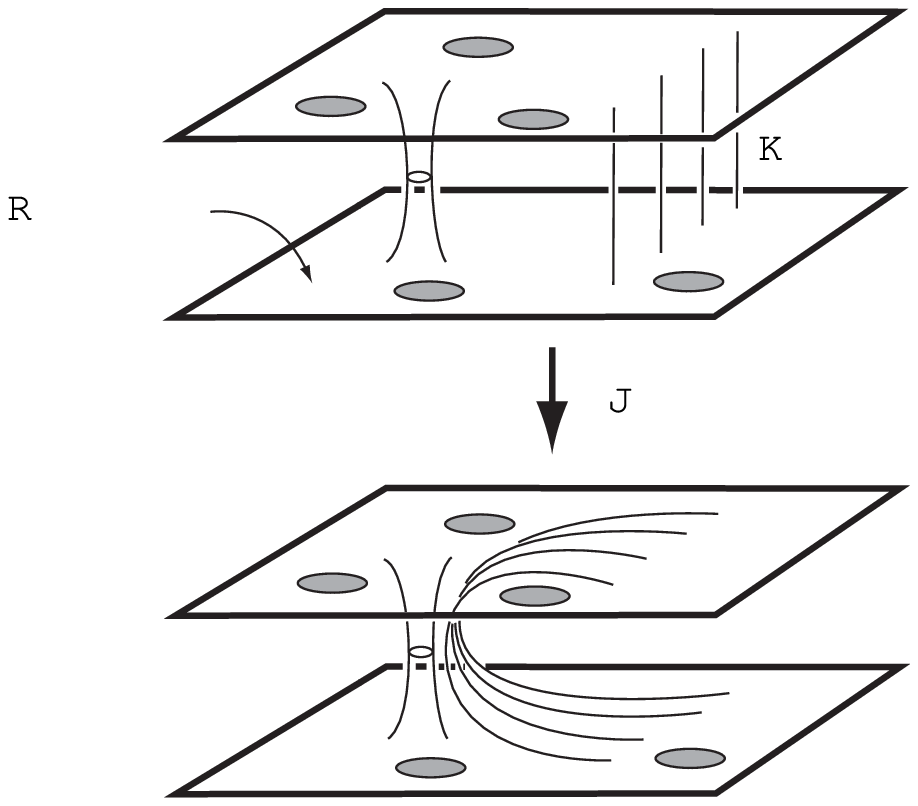}
        \caption{}
        \label{f:LightbulbIsotopy}
        \end{center}
        \end{figure}

The usual {\em lightbulb trick} (see \cite{rolfsen:90}, for example) implies that there is an isotopy, $f: M \times I \rightarrow M$, fixing the complement of $N$ pointwise such that $f_1(\alpha)$ is a straight arc (where $f_t(x)$ is shorthand for $f(x,t)$). Let $U = V \cap S$ ($ = V \setminus (\bigcup_{i=0}^{n-1} B_i \cup C_i)$). Let $K'=K \cap N$. Let $g: M \times I \rightarrow M$ be an isotopy which fixes the complement of $N$ pointwise, such that $g_1(K') \subset f_1(U)$. Such an isotopy exists since the arcs of $K'$ are straight, and since $f_1(U)$ follows the boundary of a neighborhood of $f_1(\alpha)$, which is also straight (see Figure \ref{f:LightbulbIsotopy}.) Note that $f_1^{-1}\big(g_1(K')\big) \subset U \subset S$. As before, we may now do a further isotopy of $K'$ in $S$ to bring $K'$ into $S \cap L$. 
\end{proof}

\begin{lem}
\label{l:thinballboundary}
Let $L$ be a leaf of $\mathcal F$ and $S$ be a disk properly embedded in $M$ such that 
\begin{enumerate}
	\item $\partial S$ bounds a disk on $\partial _0 M$ which contains the vertex of $K$, 
	\item $S \cap K \subset L$, and 
	\item $S \cap L$ is a connected surface.
\end{enumerate}
Let $B$ be the ball in $M$ bounded by the disk $S$ and some subdisk of $\partial _0 M$. Then either 
\begin{itemize}
	\item there exists a $K' \subset K$ such that $\partial K' \subset S \cap L$ and such that there is an isotopy of $K$ which fixes $K \backslash K'$, and sends $K'$ into $S \cap L$ or 
	\item all of the arcs of $K \cap B$ which do not lie on $\partial _0 M$ are simultaneously parallel to arcs of $K \cap B$ which do lie on $\partial _0 M$.
\end{itemize}
\end{lem}

\begin{proof}
Compress and $\partial _0$-compress $S$ as much as possible in $M^K$ to obtain a collection of disks and spheres, $R$, in $M$. If any component of $R$ is a sphere which bounds a ball containing a single unknotted arc of $K$, then there is an innermost such, $R'$. Inside $R'$ there is a disk, $T$, such that $\partial T=K' \cup \theta$, where $K'$ is a subarc of $K$ and $\theta \subset R'$. We now reverse the compressions used to obtain $R$ from $S$. Each time a compression is reversed, we attach a tube to some components of $R$. These tubes may intersect $T$, but only in its interior. We can now use $T$ to guide an isotopy of $K$ which sends $K'$ into $S$. As $\partial K' \subset S \cap L$ and $S \cap L$ is a connected subsurface of a disk we may do a further isotopy of $K'$ in $S$ to bring $K'$ into $S \cap L$. 

Since, if there are any sphere components of $R$ we are done, we now assume that only $\partial$-compressions were done to obtain $R$ from $S$. As $K$ is locally untangled we conclude that $R$ is made up of possibly several trivial disks ({\it i.e.} disks which cobound a ball, with a disk in $\partial _0 M$, containing an unknotted arc of $K$) and at least one copy of a vertex-linking disk ({\it i.e.} a disk which is parallel to the frontier of a neighborhood of the vertex of $K$). Note that there must be at least one vertex-linking disk component of $R$. (Otherwise $S$ would meet $K$ only in points on the boundary of $M$. This, in turn, implies that all of the edges of $K$ are on $\partial M$, which is not the case.) 

Let $R''$ denote a vertex-linking component of $R$ which is closest to the vertex. Every component of $K \cap B$ whose interior is in the interior of $M$ either now lies inside $R''$, or between two vertex-linking components of $R$. In either case each such subarc is  parallel to a subarc of $K$ on $\partial _0 M$.  Furthermore, undoing all of the $\partial$-compressions used to obtain $R$ from $S$ does not disturb this parallelism. 
\end{proof}

We now proceed with the proof of Theorem 6.3. Choose a locally mini-Lmax foliation $\mathcal F$ with the same minimal and maximal leaves as $\mathcal F_*$ so that when $K$ is isotoped to be mini-Lmax with respect to $\mathcal F$ the complexity $Lmax_K$ is minimal. Let $P$ be a leaf of $\mathcal F$ such that $P^K$ is a maximal leaf of $\mathcal F^K$. If $P^K$ is the first maximal leaf then let $W$ be the region of $M$ below $P$. Otherwise, let $Q$ be a leaf of $\mathcal F$ such that $Q^K$ is the minimal leaf of $\mathcal F^K$ which comes just before $P^K$, and let $W$ be the region of $M$ between $Q$ and $P$. Similarly, if $P^K$ is the last maximal leaf then let $W_*$ be the region of $M$ above $P$. Otherwise, let $Q_*$ be a leaf of $\mathcal F$ such that $Q_*^K$ is the minimal leaf of $\mathcal F^K$ which comes just after $P^K$, and let $W_*$ be the region of $M$ between $P$ and $Q_*$. 

\begin{clm}
\label{c:irreducible}
The submanifolds $W$ and $W_*$ are irreducible. 
\end{clm}

\begin{proof}
The only situation in which $W$ is not a $\partial$-compression body (and hence not irreducible) is when some component of $Q$ is a sphere. We now invoke Lemma \ref{l:thinball} (where $S=Q$) to isotope some subset of $K$ into $Q$, while keeping the rest fixed. Pushing the interior of this subset just out of $W$ shows that $Q$ can not be a minimal leaf of $\mathcal F$. 
\end{proof}

By way of contradiction, assume that there are compressing or honest $\partial _0$ compressing disks $D \subset W^K$ and $D_* \subset W^K_*$ for $P^K$ such that $D \cap D^*=\emptyset$. Note that $W^K$ is a $\partial$-compression body and $\partial _0 W^K \subset \partial _0 M \cup \partial N(K)$. 

\begin{dfn}
\rm An honest $\partial_0$-compressing disk for $P^K$ in $W^K$ or $W_*^K$ is a {\it low disk} or {\it high disk} (respectively) if $\partial D=\alpha \cup \beta$, where $\beta \subset \partial N(K)$.
\end{dfn}

\begin{dfn}
\rm A compressing or honest $\partial_0$-compressing disk for $P^K$ in $W^K$ or $W_*^K$ is {\it real} if it is a compressing or $\partial_0$-compressing disk for $P$. Otherwise it is {\it fake}.
\end{dfn}

Given this definition, we may now list all of the cases that will have to be considered (up to symmetry). 

\begin{list}{}{}
    \item (1) $D$ and $D_*$ are real.
    \item (2) $D$ and $D_*$ are fake.
    \item (3) $D$ is fake and $D_*$ is real.
\end{list}

\bigskip

\underline{Case 1.} {\it The disks $D$ and $D_*$ are real.} 

Then not only is $P^K$ a maximal leaf for $\mathcal F^K$, but also $P$ is a maximal leaf of $\mathcal F$. If $\partial D \cap \partial D_* =\emptyset$ then $P$ fails to be a strongly irreducible $\partial$-Heegaard surface, and hence, $\mathcal F$ is not locally mini-Lmax. Since the local mini-Lmaximality of $\mathcal F$ was a hypothesis of the Theorem we have reached a contradiction. This concludes Case 1.

\medskip

Following Case 1 we surmise that at least one of $D$ and $D_*$ is fake. Hence, before proceeding to the remaining cases we are forced to analyze all types of fake disks. Assume $D$ is fake. 
\begin{enumerate}
    \item $D$ is a compressing disk for $P^K$. Then $\partial D=\alpha$ is an essential curve on $P^K$, but is inessential on $P$. Hence, $\alpha$ bounds a disk $E \subset P$ which is punctured by $K$. As $W$ is irreducible (by Claim \ref{c:irreducible}), $D \cup E$ bounds a ball $B \subset W$.
        \begin{enumerate}
            \item $B$ contains the vertex of $K$. Note that this can only happen when $M$ is closed. We will call such a disk a {\it vertex compression}.
            \item $B$ does not contain the vertex. Then $B$ contains an arc of $K$ which is parallel into $P$. Hence, there is a low disk in $B$. In this case the fake disk $D$ will be referred to as a {\it cap}. A low disk in $B$ will be referred to as being {\it inside} the cap $D$. 
        \end{enumerate}
    \item $D$ is an honest $\partial _0$-compressing disk for $P^K$. Then $\partial D=\alpha \cup \beta$, where $\alpha \subset P^K$  is essential on $P^K$ and $\beta \subset \partial _0 W^K$.  
        \begin{enumerate}
            \item $\beta \subset \partial N(K)$. Then $D$ is a low disk. 
            \item $\beta \subset \partial _0 M$. Then there is a disk $E \subset P$ punctured by $K$ such that $\partial E=\alpha \cup \gamma$, where $\gamma \subset \partial P$. If $\beta$ was essential on $\partial _0 W$ then $D \cup E$ would be a compressing disk for $\partial _0 W$, which is impossible. (To see the contradiction, double $W$ along $\partial _0 W$. The compressing disk $D \cup E$ then turns into an essential 2-sphere in a compression body.) We conclude $\beta$ is inessential on $\partial _0 W$, and hence cobounds a disk $V$. If $V \cap K =\emptyset$ then $D$ would not be an honest $\partial _0$-compressing disk for $P^K$. Hence, $V$ either contains a low disk or it contains the vertex of $K$. In the former case the disk $D$ will be referred to as a {\it $\partial$-cap}. In the latter it will be referred to as a {\it vertex $\partial$-compression}. If $D$ is a $\partial$-cap then a low disk in $V$ is said to be {\it inside} $D$. 
        \end{enumerate}
\end{enumerate}

\medskip
\newpage

\underline{Case 2.} {\it $D$ and $D_*$ are fake.}

\medskip

This case is organized into subcases as follows:

\begin{list}{}{}
        \item (2.1) $D$ and $D_*$ are compressing disks. 
        \item (2.2) $D_*$ is a compressing disk and $D$ is a $\partial$-compressing disk. 
        \item (2.3) $D$ and  $D_*$ are $\partial$-compressing disks.
\end{list}

\medskip

{\it Subcase 2.1.} If $D$ and $D_*$ are compressing disks then $\partial D$ bounds a disk $E$ in $P$ and $\partial D_*$ bounds a disk $E_*$ in $P$. As $D \cap D_* =\emptyset$ the following is a complete list of the subcases that will have to be considered, up to symmetry.

\begin{list}{}{}
    \item (2.1.1) $E \cap E_*=\emptyset$
    \item (2.1.2) $E_* \subset E$
\end{list}

        \begin{figure}[htbp]
        \vspace{0 in}
        \begin{center}
        \includegraphics[width=3 in]{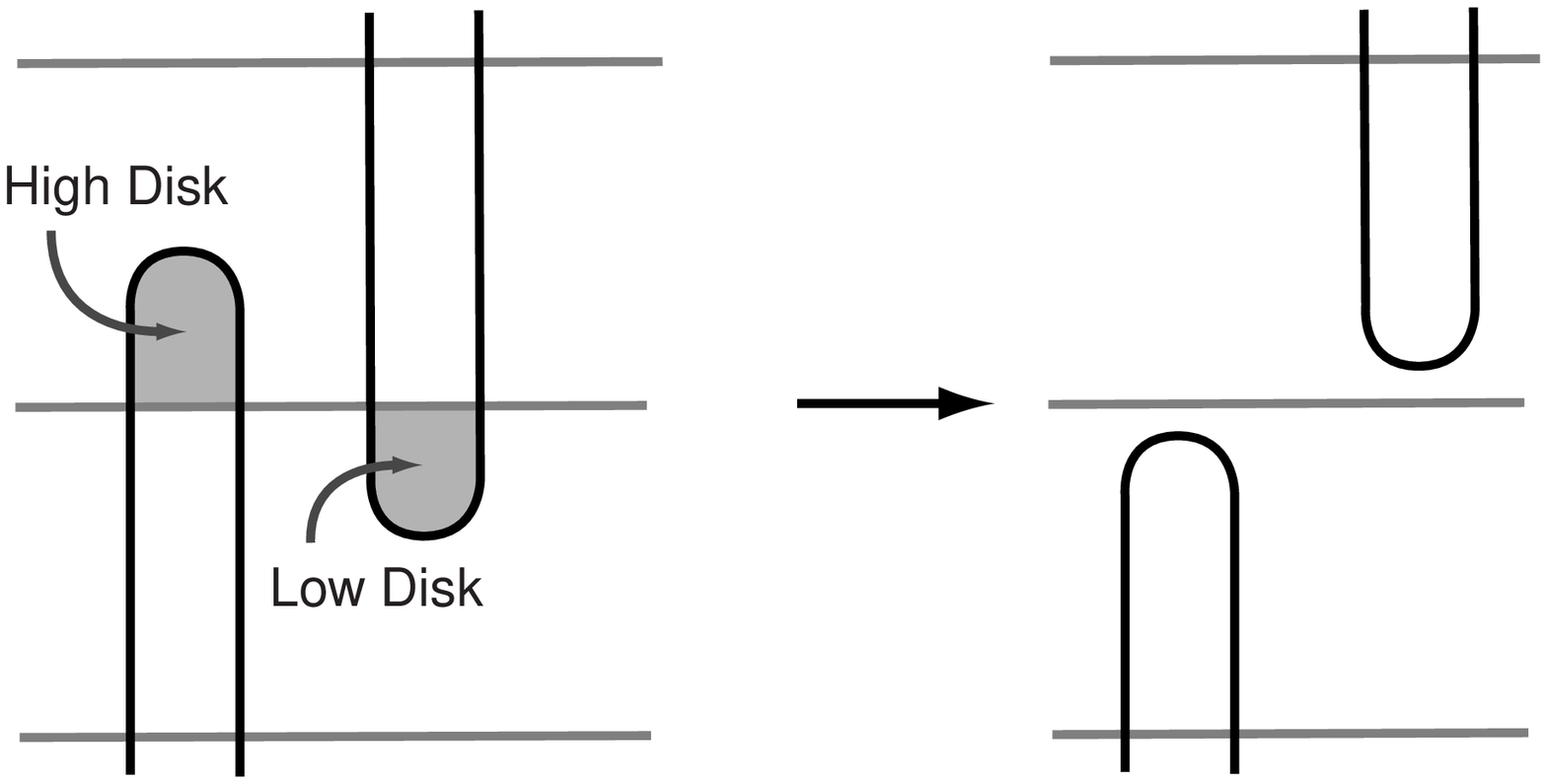}
        \caption{}
        \label{f:High&Low}
        \end{center}
        \end{figure}

\begin{figure}[htbp]
        \vspace{0 in}
        \begin{center}
        \includegraphics[width=3 in]{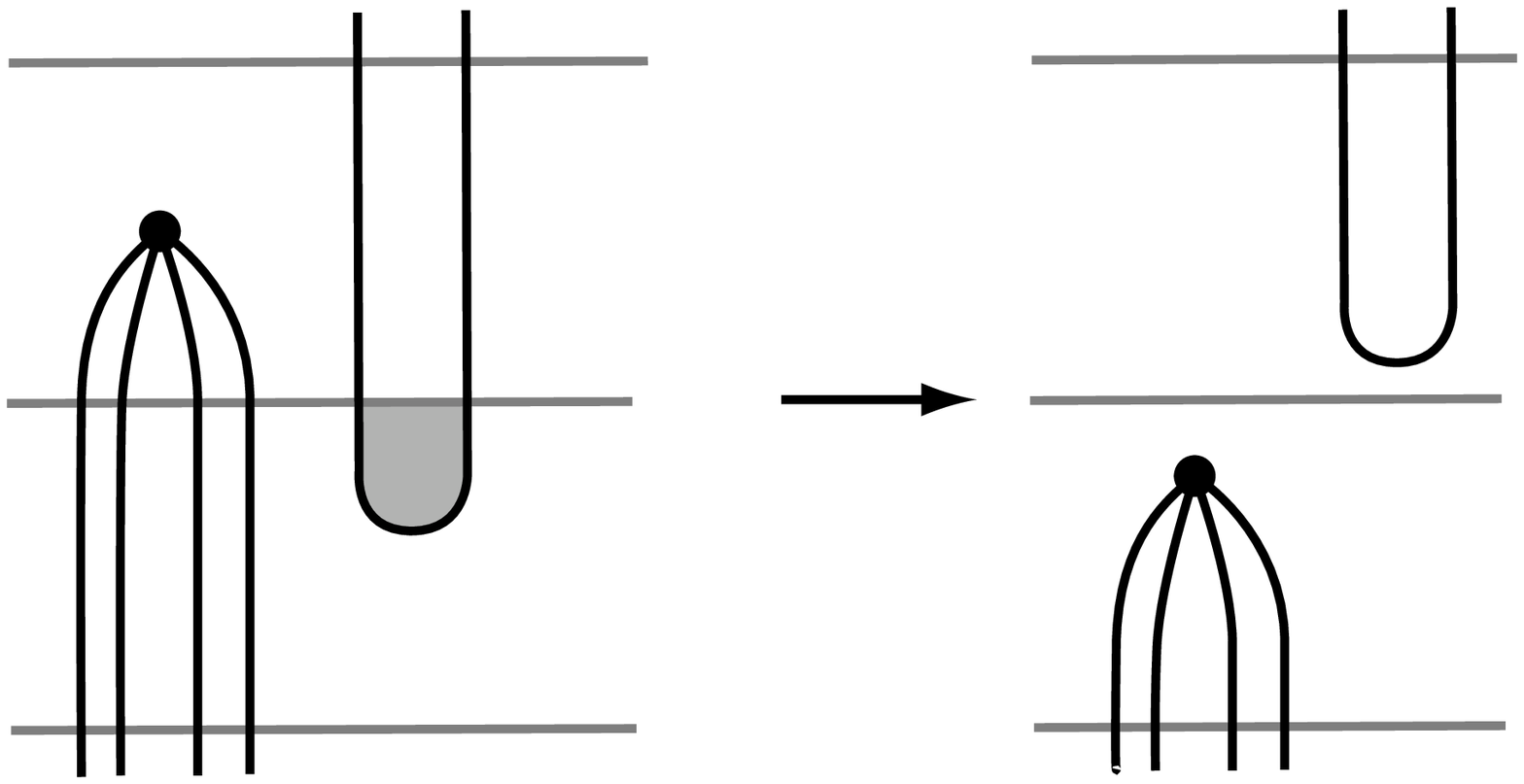}
        \caption{}
        \label{f:Low+Vertex}
        \end{center}
        \end{figure}

\medskip
{\it Subcase 2.1.1.} $E \cap E_*=\emptyset$. If $D$ and $D_*$ are caps then there is a low disk inside $D$ which is disjoint from a high disk inside $D_*$. We can then do the move depicted in Figure \ref{f:High&Low} to reduce $Lmax_K(\mathcal F)$. On the other hand, if one of $D$ or $D_*$ is a vertex compression then we can perform the move depicted in Figure \ref{f:Low+Vertex} to lower $Lmax_K(\mathcal F)$. 

\medskip
{\it Subcase 2.1.2.} $E_* \subset E$. Then $\partial D$ and $\partial D_*$ cobound an annulus $A \subset P$. We assume that $D_*$ is outermost, in the sense that there does not exist a compressing disk $D_*' \subset W_*$ for $P^K$ such that $\partial D_*' \subset A$, but is not parallel to $\partial D_*$ in $P^K$. Let $S$ denote the sphere $D \cup A \cup D_*$. We now claim that $S^K$ is incompressible in $M^K$. Suppose not.

Invoke Lemma \ref{l:thinball} to isotope some subset $K'$ of $K$ into $A$. If the interior of $K'$ met at least one maximal leaf of $\mathcal F^K$ then this isotopy  produces a foliation with smaller $Lmax_K$. Furthermore, since we have only done an isotopy of $K$ we have not disturbed the foliation $\mathcal F$, so we have reached a contradiction. 

If the interior of $K'$ does not meet any maximal leaves of $\mathcal F^K$ then it also must not meet any minimal leaves. Hence, every arc of $K'$ lies in $W$ or $W_*$. If some such arc, $\gamma$, were in $W$ then it would also be in the ball bounded by $D \cup E$. As it is isotopic to an arc in $A$ there is a low disk $L \subset W$ such $L \cap P \subset A$. But then this low disk is disjoint from $E_*$, and we may do one of the moves depicted in Figure \ref{f:High&Low} or \ref{f:Low+Vertex}, according to whether $D_*$ is a cap or a vertex compression. 

We conclude $\gamma \subset W_*$. As it is isotopic to an arc in $A$ there is a high disk $H \subset W_*$ such $H \cap P \subset A$. Let $H'$ be a cap in $W_*$ whose boundary is in $A$, such that the high disk $H$ is inside $H'$. Let $s$ be any arc in $A^K$ joining $\partial H'$ to $\partial D_*$. Then there is a component of the frontier of a neighborhood of $H' \cup s \cup D_*$ in $W_*$ which is a cap which is ``more outermost" than $D_*$, contradicting our assumption. 

        \begin{figure}[htbp]
        \psfrag{D}{$D$}
        \psfrag{d}{$D_*$}
        \vspace{0 in}
        \begin{center}
        \includegraphics[width=2 in]{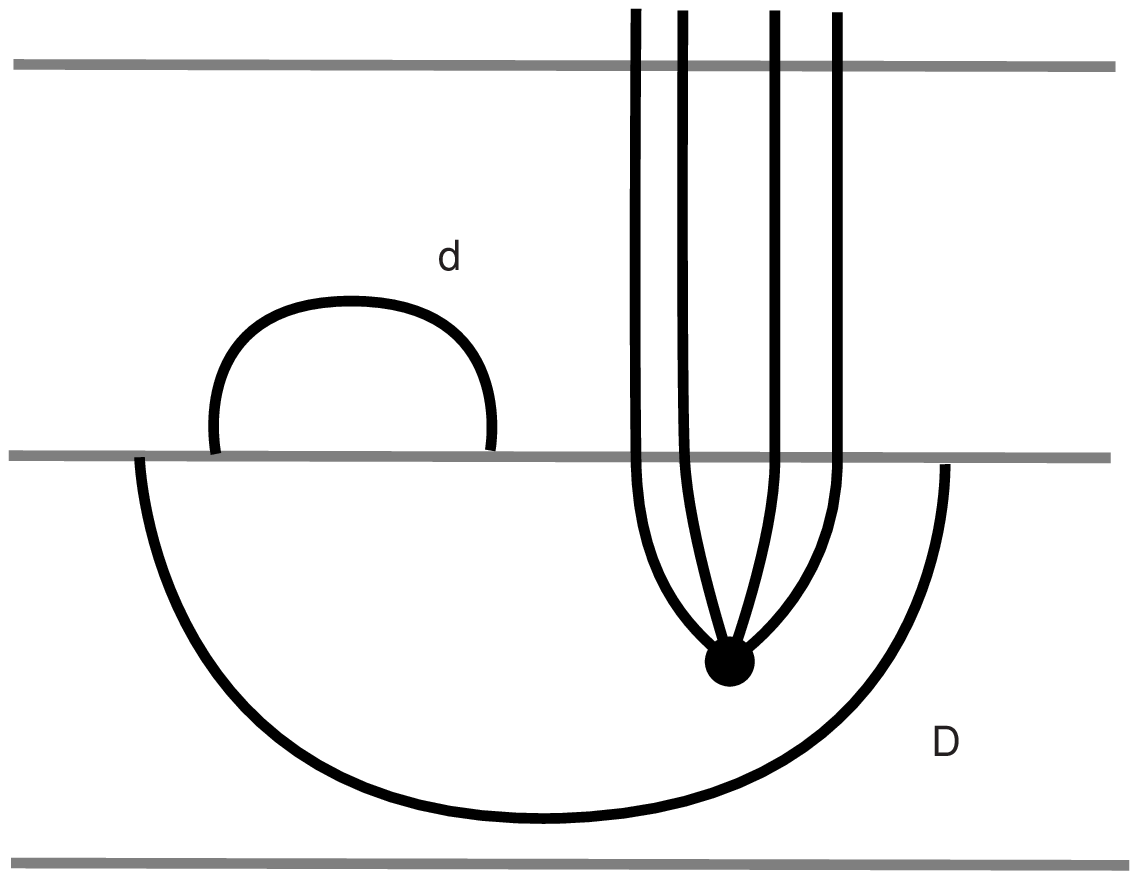}
        \caption{}
        \label{f:VertexOnly}
        \end{center}
        \end{figure}

We are now faced with the possibility that $S^K$ was incompressible to begin with. As the only such spheres are parallel to the vertex of $K$ we may now isotope $K$ to look like Figure \ref{f:VertexOnly} in the ball bounded by $S$. As this isotopy necessarily reduced the number of intersections of $K$ with $P$, the complexity $Lmax_K$ has gone down.

\medskip
{\it Subcase 2.2.} If $D_*$ is a compressing disk and $D$ is a $\partial$-compressing disk then $\partial D_*$ bounds a disk $E_* \subset P$ and $\partial D=\alpha \cup \beta$ where $\alpha \subset P$. There are now three further subcases. 

                \begin{list}{}{}
                    \item (2.2.1) $D$ is a $\partial$-cap
                    \item (2.2.2) $D$ is a low disk
                    \item (2.2.3) $D$ is a vertex $\partial$-compression
                \end{list}

\medskip
{\it Subcase 2.2.1} If $D$ is a $\partial$-cap then there is a low disk $D'$ inside $D$ which is disjoint from $E_*$. If $D_*$ is a cap then there is a high disk inside $D_*$, which is then disjoint from $D'$. We can now do the $Lmax$ lowering move depicted in Figure \ref{f:High&Low}. Similarly, if $D_*$ is a vertex compression then we can do the move depicted in Figure \ref{f:Low+Vertex}. 

\medskip
{\it Subcase 2.2.2} If $D$ is a low disk then there is a cap $D'$ such that $D$ is inside $D'$ and $D' \cap D=\emptyset$. This case now reduces to Subcase 2.1. 

\medskip
{\it Subcase 2.2.3} The disk $D$ is a vertex $\partial$-compression, and hence $\beta \subset \partial _0 W$. It follows that $D_*$ is a cap, since the vertex of $K$ is in $W$.  The curves $\partial D$ and $\partial D_*$ cut off subdisks $E$ and $E_*$ of $P$. 

\begin{figure}[htbp]
        \psfrag{D}{$D_*$}
        \psfrag{d}{$D$}
        \psfrag{f}{$D'$}
        \vspace{0 in}
        \begin{center}
        \includegraphics[width=5.5 in]{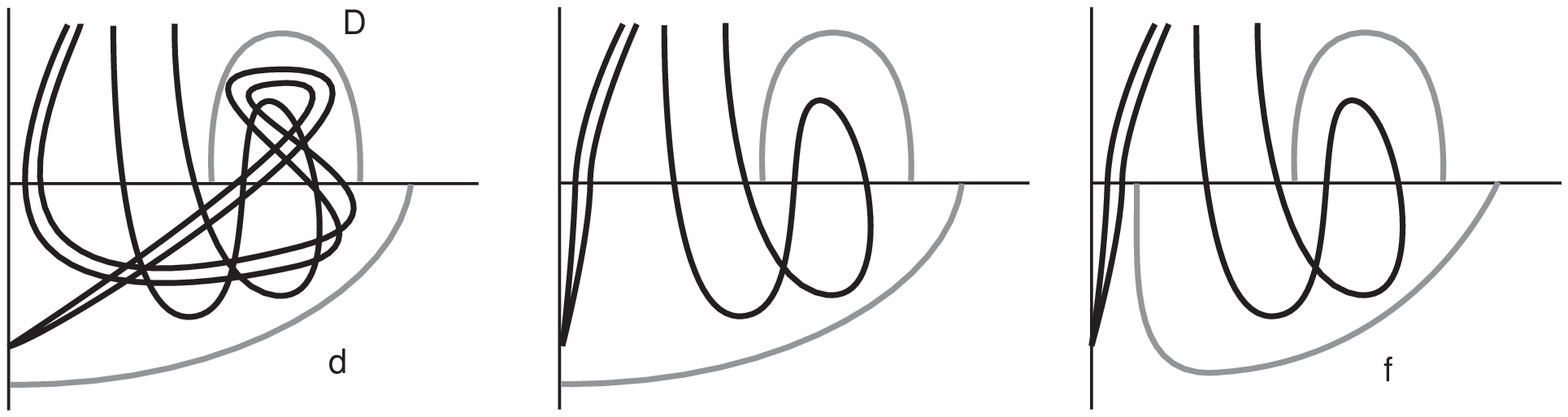}
        \caption{}
        \label{f:boundarycompression/compression}
        \end{center}
        \end{figure}

If $E \cap E_* =\emptyset$ then there is a high disk inside $D_*$ which we can use to perform the $Lmax_K$ lowering move depicted in Figure \ref{f:Low+Vertex}. Hence, $E_* \subset E$, and we are in a situation similar to Subcase 2.1.2 (see Figure \ref{f:boundarycompression/compression} left). Let $A$ be the closure of $E-E_*$. We assume that $D_*$ is outermost, in the sense that there does not exist a compressing disk $D_*' \subset W_*$ for $P^K$ such that $\partial D_*' \subset A$, but is not parallel to $\partial D_*$ in $P^K$. Let $S$ denote the disk $D \cup A \cup D_*$. 

By Lemma \ref{l:thinballboundary} there are two possibilities. The first is that some subarc $K'$ of $K$ is isotopic into $A$. If the interior of $K'$ met at least one maximal leaf of $\mathcal F^K$ then this isotopy produces a foliation with smaller $Lmax_K$. Furthermore, since we have only done an isotopy of $K$ we have not disturbed the foliation $\mathcal F$, so we have reached a contradiction. 

If the interior of $K'$ does not meet any maximal leaves of $\mathcal F^K$ then it also must not meet any minimal leaves. Hence, $K'$ lies in $W$ or $W_*$. If $K'$ were in $W$ then it would also be in the ball bounded by $D \cup E \cup V$. As it is isotopic to an arc in $A$ there is a low disk $L \subset W$ such $L \cap P \subset A$. But then this low disk is disjoint from $E_*$, and we may do the move depicted in Figure \ref{f:High&Low}. 

We conclude $K' \subset W_*$. As it is isotopic to an arc in $A$ there is a high disk $H \subset W_*$ such $H \cap P \subset A$. Let $H'$ be a cap in $W_*$ whose boundary is in $A$, such that the high disk $H$ is inside $H'$. Let $s$ be any arc in $A^K$ joining $\partial H'$ to $\partial D_*$. Then there is a component of the frontier of a neighborhood of $H' \cup s \cup D_*$ in $W_*$ which is a cap which is ``more outermost" than $D_*$, contradicting our assumption. 

The second possibility implied by the conclusion of Lemma \ref{l:thinballboundary} is that the subarcs $K''$ of $K$ which join the vertex of $K$ to $A$, and whose interiors lie in the interior of $M$, are isotopic to subarcs of $K$ on $\partial _0 M$  (see Figure \ref{f:boundarycompression/compression} middle). If any of the arcs of $K''$ met $E_*$ then this isotopy lowers the complexity $Lmax_K$. If not then we see there is a fake compressing disk, $D'$, whose boundary is disjoint from $D_*$ (see Figure \ref{f:boundarycompression/compression} right), and we are reduced to Subcase 2.1.2.

\medskip
{\it Subcase 2.3.} If $D$ and $D_*$ are $\partial$-compressing disks then the following table lists all subcases, up to symmetry. 

\medskip

\begin{center}
\begin{tabular}{l|l|l|l}
$D \backslash D_*$ & $\partial$-cap & high disk & vertex $\partial$-compression\\
\hline
$\partial$-cap & (2.3.1) & (2.3.1) & (2.3.3) \\
\hline
low disk &  &(2.3.1) & (2.3.2)\\
\hline
vertex $\partial$-compression & & & Not possible
\end{tabular}
\end{center}

\medskip

\begin{figure}[htbp]
        \vspace{0 in}
        \begin{center}
        \includegraphics[width=3 in]{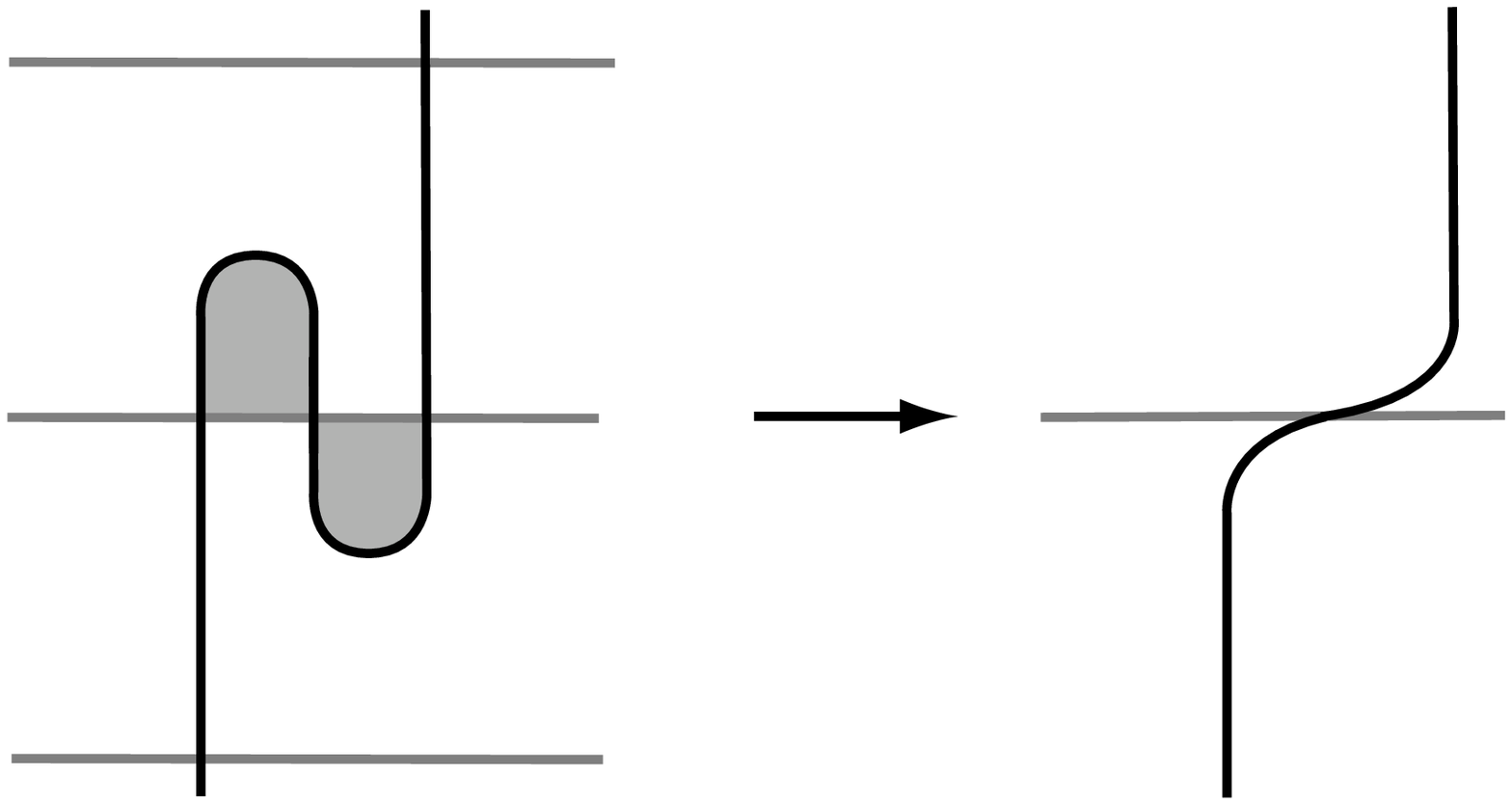}
        \caption{}
        \label{f:SCurve}
        \end{center}
        \end{figure}

{\it Subcase 2.3.1.} If $D$ is a low disk or a $\partial$-cap and $D_*$ is a high disk or a $\partial$-cap then there is a low and a high disk that we can use to perform one of the moves depicted in Figures \ref{f:High&Low} or \ref{f:SCurve}. 

\medskip

{\it Subcase 2.3.2.} If $D$ is a low disk in the interior of $M$ and $D_*$ is a vertex $\partial$-compression then there is a cap $D'$ in $W$ such that $D$ is inside $D'$, and $\partial D' \cap \partial D_*=\emptyset$. We are now reduced to Subcase 2.2.3. If $D \subset \partial _0 M$ then there is a $\partial$-cap $D'$ which contains $D$ such that $\partial D' \cap \partial D_*=\emptyset$. 
This is handled in Subcase 2.3.4.

\medskip

{\it Subcase 2.3.3.} The final case is when $D_*$ is a vertex $\partial$-compression and $D$ is a $\partial$-cap. Then $\partial D$ and $\partial D_*$ cut off subdisks $E$ and $E_*$ of $P$. Inside $D$ there is a low disk. If $E \cap E_*=\emptyset$ then we may use this low disk to perform the move depicted in Figure \ref{f:Low+Vertex}. We conclude $E_* \subset E$ or $E \subset E_*$. Let $S$ denote the disk comprised of $E$, $E_*$ and the subset of $P$ which lies between $E$ and $E_*$. 

Just as in Subcase 2.2.3 we now invoke Lemma \ref{l:thinballboundary} with the disk $S$. As before the first possibility is that some subarc $K'$ of $K$ is isotopic into $A$. We rule this possibility out by reasoning which is identical to that in Subcase 2.2.3.

\begin{figure}[htbp]
        \psfrag{D}{$D$}
        \psfrag{d}{$D_*$}
        \vspace{0 in}
        \begin{center}
        \includegraphics[width=4 in]{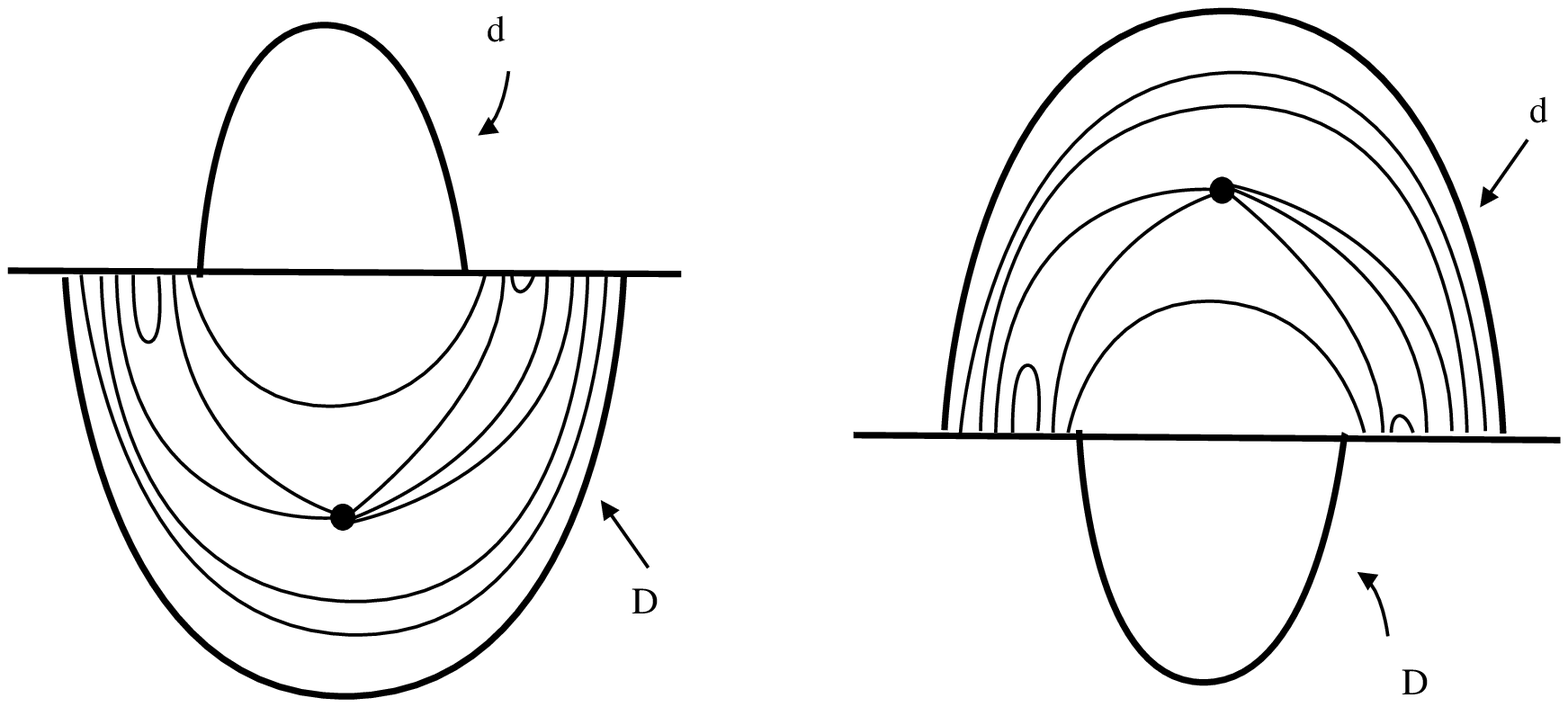}
        \caption{}
        \label{f:C233}
        \end{center}
        \end{figure}

Let $B$ denote the ball bounded by $S$ and a subdisk of $\partial _0 M$. The second possibility given by Lemma \ref{l:thinballboundary} is that all of the arcs $K''$ of $K \cap B$ which do not lie on $\partial _0 M$ are isotopic to subarcs of $K$ which do lie on $\partial _0 M$. But we can always isotope the arcs on $\partial _0 M$ to look like one of the two possibilities given in Figure \ref{f:C233} (depending on whether $E \subset E_*$ or $E_* \subset E$). Dragging the arcs $K''$ along with such an isotopy reduces $Lmax _K$. Note that in the figure on the left the vertex of $K$ has crossed $P$. 

\medskip

\underline{Case 3.} {\it $D$ is fake and $D_*$ is real.}

\medskip

In $W_*$ we alter the foliation $\mathcal F$ as follows. Note that $Q_*$ is obtained from $P$ by a sequence of compressions. These compressions are realized by the foliation $\mathcal F$ in the sense that as $t$ (the height parameter for $\mathcal F$) increases we see the leaves compress in some sequence. But note that any sequence of compressions that one can perform on $P$ to obtain $Q_*$ can be realized by some foliation in this way. Choose some such sequence which begins with a compression along the disk $D_*$. Alter the foliation in $W_*$ so that the leaves compress along this sequence. Note that the maximal and minimal leaves of $\mathcal F^K$ are unchanged, so the complexity $Lmax_K$ is unchanged as well.

The disk $D_*$ now represents a compression or $\partial_0$-compression of the leaves of $\mathcal F$. If $D$ is a low disk then we may use it to isotope a subarc of $K$ so that we pass through the corresponding minimum of $K$ {\it after} we see the compression/$\partial _0$-compression which corresponds to $D_*$. Such a move reduces the complexity $Lmax_K$, a contradiction.

If $D$ is a compression then we let $B$ be the ball bounded by $D$ and some subdisk of $P$. If $D$ is a $\partial$-compression then let $B$ be the ball bounded by $D$, a subdisk of $P$, and a subdisk of $\partial _0 W$. In either case we can use $B$ to define an isotopy of $K$ so that any minimum or vertex of $K$ which is inside $B$ gets pushed past the compression/$\partial _0$-compression which corresponds to $D_*$. Again, such a move reduces the complexity $Lmax_K$.

\medskip

In summary, we have shown that if $D$ is any compressing (or honest $\partial _0$-compressing) disk for $P^K$ in $W^K$, and $D_*$ is a compressing (or honest $\partial _0$-compressing) disk for $P^K$ in $W_*^K$, then $\partial D \cap \partial D_* \ne \emptyset$. Hence, $P^K$ is a quasi-strongly irreducible $\partial$-Heegaard surface for $(W^K) \cup (W_*^K)$. Furthermore, the only times the assumption of honesty was used was in cases where $K$ had a vertex. Hence, if $K$ has no vertex then we conclude $P^K$ is strongly irreducible.

\bibliographystyle{alpha}